# MULTIPROCESS PARALLEL ANTITHETIC COUPLING FOR BACKWARD AND FORWARD MARKOV CHAIN MONTE CARLO[1]

By Radu V. Craiu and Xiao-Li Meng

*University of Toronto and Harvard University*

Antithetic coupling is a general stratification strategy for reducing Monte Carlo variance without increasing the simulation size. The use of the antithetic principle in the Monte Carlo literature typically employs two strata via antithetic quantile coupling. We demonstrate here that further stratification, obtained by using $k > 2$ (e.g., $k = 3$–10) antithetically coupled variates, can offer substantial additional gain in Monte Carlo efficiency, in terms of both variance and bias. The reason for reduced bias is that antithetically coupled chains can provide a more dispersed search of the state space than multiple independent chains. The emerging area of perfect simulation provides a perfect setting for implementing the $k$-process parallel antithetic coupling for MCMC because, without antithetic coupling, this class of methods delivers genuine independent draws. Furthermore, antithetic backward coupling provides a very convenient theoretical tool for investigating antithetic forward coupling. However, the generation of $k > 2$ antithetic variates that are *negatively associated*, that is, they preserve negative correlation under monotone transformations, and *extremely antithetic*, that is, they are as negatively correlated as possible, is more complicated compared to the case with $k = 2$. In this paper, we establish a theoretical framework for investigating such issues. Among the generating methods that we compare, Latin hypercube sampling and its iterative extension appear to be general-purpose choices, making another direct link between Monte Carlo and quasi Monte Carlo.

**1. Paired antithetic coupling.** Monte Carlo estimation of the expectation of an estimand function $f$ with respect to a probability measure $\pi$, $\mu = \int f(x)\pi(dx)$, can be done in many ways. The simplest method, known

Received October 2002; revised November 2003.
[1]Supported in part by NSF, NSA and NSERC grants.

as "crude Monte Carlo," proceeds by simulating a (not necessarily independent) sample $X_1, X_2, \ldots, X_n$ from $\pi$ and estimating $\mu$ by the sample average $\hat{\mu}_n = \frac{1}{n} \sum_i f(X_i)$. More refined methods, collectively known as *swindles*, take advantage of well-known statistical principles to construct more efficient designs and/or more efficient estimators. In this paper we focus on antithetic coupling, which can be viewed as a way to induce stratification. For other type of swindles see Hammersley and Handscomb ([1965](#)), Simon ([1976](#)), Kennedy and Gentle ([1980](#)) and Gentle ([1998](#)), among others.

In the context of classic independent sampling, antithetic coupling [Hammersley and Morton ([1956](#))] is commonly described as a method of producing two negatively correlated copies of an unbiased estimator. The average of the two, each based on $m$ draws, will then have smaller variance compared to the estimator based on the $n = 2m$ independent draws. Substantial reductions have been documented throughout the literature from early works, such as Page ([1965](#)) and Fishman ([1972](#)), to most recent ones such as Frigessi, Gåsemyr and Rue ([2000](#)), who demonstrate the power of the antithetic principle in the more complicated (forward) MCMC setting.

The negative correlation between the two copies is typically induced via antithetic quantile coupling by using a pair $\{U, 1 - U\}$, where $U \sim$ Uniform$(0, 1)$, in the sample generating process. The amount of variance reduction is governed by the degree of symmetry in the distribution of our estimator. In the extreme case when this distribution is symmetric, the use of paired antithetic variates can entirely eliminate the Monte Carlo variance when the underlying draws are independent or reduce the variance from the usual $n^{-1}$ rate to $n^{-2}$ rate when the draws are realizations of a genuine MCMC algorithm, as observed by Frigessi, Gåsemyr and Rue ([2000](#)). This emphasizes the stratification aspect of the antithetic principle as a way to divide the sample space into a "negative" stratum and a "positive" stratum. The equal amount of draws from each stratum ensured through *pairing* brings in further variance reduction compared to using simple random sampling within each stratum. Generalizing the antithetic coupling from $k = 2$ to $k > 2$ processes is quite natural from this stratification point of view, as often more than two strata can produce substantial additional gain. A main goal of this paper is to demonstrate such gains in the context of MCMC, as well as the additional benefit of improving the mixing of the original Markov chains.

The paired quantile coupling via $\{U, 1 - U\}$ has the following *extreme antithesis* (EA) property that is usually not emphasized in the literature. Specifically, if $F$ is an *arbitrary* univariate cumulative distribution function (CDF), and $X_1 = F^{-1}(U)$, $X_2 = F^{-1}(1 - U)$, where $U \sim$ Uniform$(0, 1)$, then Corr$(X_1, X_2)$ achieves *the minimal possible value* subject to the constraint that $X_1, X_2 \sim F$. The proof [Moran ([1967](#))] relies on the elegant Hoeffding

identity

$$\text{(1.1)} \qquad \text{Cov}(X_1, X_2) = \iint [F(x_1, x_2) - F(x_1)F(x_2)]\mu(dx_1\, dx_2),$$

and the equally elegant Fréchet (1951) inequality

$$\text{(1.2)} \quad \max\{F(x_1) + F(x_2) - 1, 0\} \le F(x_1, x_2) \le \min\{F(x_1), F(x_2)\},$$

where $F(x_1, x_2)$ is the joint CDF. The fact that the single strategy of using $\{U, 1-U\}$ achieves EA *simultaneously for all* $F$ also implies that it achieves EA for $\text{Corr}(g(X_1), g(X_2))$, where $g$ is any monotone function such that $\int g^2(x)F(dx) < \infty$. It is known that negative correlation is not even preserved by monotone transformations, so the above discussion implies that $\{U, 1-U\}$ must satisfy a stronger condition, *negative association* (NA) [Joag-Dev and Proschan (1983)], which requires exactly that the negative correlation be preserved by monotone functions of the variates (see Section 3.1). The fact that $\{U, 1-U\}$ is automatically NA appears to be responsible for the general silence of the NA requirement in the literature of antithetic coupling, but once we move beyond $k = 2$, NA (as well as EA) becomes a key notion in our general theoretical foundation.

The rest of the paper is divided into three sections. Section 2 presents the general $k$-process parallel antithetic coupling technique, illustrated with a backward MCMC mixture sampling and two forward MCMC sampling schemes. However, for $k > 2$, as we demonstrate in Section 3, there is no general strategy that achieves EA simultaneously even just for uniform and normal distributions. As a result, it is harder to ensure NA and EA for a given problem, especially in theory. Section 3 is thus devoted to a general theory for ensuring NA and EA for arbitrary $k$, which is our main theoretical contribution. Section 4 discusses several common methods for generating $k$ antithetic variates and their general properties with respect to achieving EA and NA.

## 2. $k$-process parallel antithetic coupling and illustration.

2.1. *Perfect simulation and time-backward dual sequence.* Since the seminal work by Propp and Wilson (1996) on *coupling from the past* (CFTP), there has been an array of methodological and theoretical papers on how to use backward coupling methods for exact MCMC sampling. David Wilson's web site (http://dimacs.rutgers.edu/~dbwilson/exact) is the most comprehensive single source for learning about the fast-growing field of *exact sampling* or *perfect simulation*, so named because, for this class of sampling methods, the thorny issue of deciding the running time for an acceptable error in approximating the distribution of interest, $\pi$, completely vanishes. A CFTP algorithm, or Fill's algorithm (1998), or many of their variations and generalizations [e.g., Fill, Machida, Murdoch and Rosenthal (2000), Meng (2000) and Wilson (2000b)] will terminate itself with probability 1 in a finite amount of time and yet delivers authentic (and hence *exact/perfect*) independent draws from the limiting distribution.

This can be achieved by constructing, for a Markov chain $\{X_t\}_{t \geq 0}$ with stationary distribution $\pi$, the "dual" sequence $\tilde{X}_t$ such that it has the same distribution as $X_t$ but it purposely violates the Markovian property. This is perhaps most easy to illustrate by first considering *a time-inhomogeneous* Markov chain defined by $X_t = \psi_t(X_{t-1}, U_t)$, where $\{U_t, t \in \mathcal{N}\}$ are i.i.d. random variables. For any $t > 0$ and any starting value $X^*$, we can compute a *time-forward* sequence

$$(2.1) \qquad X_t = \psi_t(\psi_{t-1}(\ldots \psi_1(X^*, U_1) \ldots, U_{t-1}), U_t),$$

as well as a *time-backward* sequence

$$(2.2) \qquad \tilde{X}_t = \psi_1(\psi_2(\ldots \psi_t(X^*, U_t) \ldots, U_2), U_1).$$

Clearly, $X_t$ is Markovian because $X_t = \psi_t(X_{t-1}, U_t)$, but $\tilde{X}_t$ is not.

Evidently, when $\psi_t$ is time-homogeneous, $X_t$ and $\tilde{X}_t$ have identical distributions for any $t$. What is gained by giving up the Markovian property is that the backward sequence $\tilde{X}_t$ can "hit and stay at" the limit [Thorisson (2000), Chapter 1] in a finite amount of time when its forward counterpart $X_t$ cannot. This is easiest to see in an extreme case where $\psi(X, U) = U$, which, say, is uniform on (0, 1). Then it is obvious that $X_t = U_t$, but $\tilde{X}_t = U_1$, for all $t \geq 1$. In other words, while both $X_t$ and $\tilde{X}_t$ converge, trivially, to Uniform(0, 1) in distribution, $X_t$ wanders off from $U_t$ to $U_{t+1}$, even if it hits the right limiting distribution already at $t = 1$. In contrast, $\tilde{X}_t$ stays at the exact same value $U_1$ for all $t \geq 1$. In general, assuming $X_t$ is uniformly ergodic, one can conclude that while $X_t$ converges *in distribution*, there exist an updating $\psi$ and a finite stopping time $T$ such that $\tilde{X}_T$ hits the same limit *with probability* 1 [Foss and Tweedie (1998)]. By mapping $t$ to $-t$ in (2.2), and thereby creating a convenient forward execution of the time-backward sequence given in (2.2), Propp and Wilson (1996) devise the CFTP method for finding such a $T$.

At the crux of CFTP implementation is the ability to follow a large, possibly infinite, number of paths in time and to assess whether these paths have all merged after a certain time $T$. For a good introduction to CFTP we refer to Casella, Lavine and Robert (2001). In real applications, the most challenging problem is that for many routine statistical problems, such as posterior sampling, the state space is both uncountable and unbounded. Although intense work has been done in this area, in terms of both general strategies and specific implementation [e.g., Kendall (1998), Murdoch and Green (1998), Møller and Nicholls (1999), Møller (1999), Murdoch (2000),

Hobert, Robert and Titterington (1999), Casella, Mengersen, Robert and Titterington (2002) and Philippe and Robert (2003)], the difficulties exhibited in Murdoch and Meng (2001) in the context of posterior sampling with a $t$ likelihood and normal mixture priors indicate that much more research is needed before CFTP and its variants and extensions can become a method of choice in common statistical applications.

2.2. *Multiprocess antithetically coupled CFTP.* Since CFTP delivers i.i.d. draws from the target density, when viewed as a "black box," it is no different from many classical Monte Carlo sampling methods, such as inverse CDF transform or rejection sampling. It is then natural to consider antithetic coupling for CFTP. However, unlike classical methods, the CFTP black box is a mapping from an infinite product space for $\{U_t, t \le 0\}$ to the state space $S$, which makes the underlying theory for guaranteeing NA (and EA) of the samples more complicated than those available in the literature. A theoretical foundation is therefore needed, and this is established in Section 3. Here we focus on the description of the method itself.

At the core of our method is the generation of $k$ negatively associated $\{U^{(1)}, \ldots, U^{(k)}\}$, each having the same marginal distribution as the $U$ required by $\psi(\cdot, U)$. There are many ways for doing so for a given distribution of $U$; see Section 4. Given such a method, one can run $k$ CFTP processes in parallel, the $j$th one using $\{U_t^{(j)}, t \le 0\}$, $j = 1, \ldots, k$, where $\{U_t^{(1)}, \ldots, U_t^{(k)}\}$, $t \le 0$, are i.i.d. copies of $\{U^{(1)}, \ldots, U^{(k)}\}$, as sketched in Figure 1. Within the $j$th process of CFTP all paths are positively coupled because they use the same $\{U_t^{(j)}, t \le 0\}$. At each update, $\{U_t^{(1)}, \ldots, U_t^{(k)}\}$ are NA, a property that clearly does not alter the validity of each individual CFTP process.

To obtain $n = km$ draws, we repeat the above procedure *independently* $m$ times and collect $\{X_i^{(j)}, 1 \le i \le m; 1 \le j \le k\}$, where $i$ indexes the replication, as our sample $\{X_1, X_2, \ldots, X_n\}$. Let $\sigma_f^2 = \text{Var}_\pi[f(X)]$ and $\rho_k^{(f)} = \text{Corr}_\pi(f(X_1^{(1)}), f(X_1^{(2)}))$, which is intended to be negative. Then

$$(2.3) \qquad \text{Var}\left(\frac{1}{n}\sum_{i=1}^{n} f(X_i)\right) = \frac{\sigma_f^2}{n}[1 + (k-1)\rho_k^{(f)}].$$

Consequently, the *variance reduction factor* (VRF) relative to the independent sampling *with the same simulation size*, is

$$(2.4) \qquad S_k^{(f)} = 1 + (k-1)\rho_k^{(f)}.$$

We emphasize here the dependence of $S_k^{(f)}$ on $k$ and more importantly on $f$, and thus the actual gain in reduction can be of practical importance for some $f$ but not for others.

The size-fixed VRF $S_k^{(f)}$ ignores the possible different computational requirements between generating $k$ independent draws and $k$ antithetically coupled draws. A completely fair comparison is typically impossible with any simulation study because the overall computational time is seldom exactly linear in size and it can depend critically on software, hardware, programming skill and the actual implementation. Nevertheless, a useful approximation can be derived by assuming linearity and by ignoring any "overhead." Specifically, let $\tau_k$ be the average CPU time needed to make a *joint* draw from a particular antithetic $k$-process. Then, under our assumptions, for a given total CPU time $T$, the average number of dependent draws we can make is $n_{\text{dep}} = Tk/\tau_k$, and the number of independent draws is $n_{\text{indep}} = T/\tau_1$. Consequently, the *time-fixed* VRF is given by

$$(2.5) \qquad T_k^{(f)} = C_k S_k^{(f)} \qquad \text{where } C_k = \frac{\tau_k/k}{\tau_1}.$$

Note $C_k$ is free of $f$, but a large $C_k$ may offset the gain in $S_k^{(f)}$ (i.e., by making $T_k^{(f)} > 1$), as seen in the next section. Also note that although we use the $\tau_1$ notation for consistency, there can be a substantial difference between setting $k = 1$ in a general $k$-process subroutine and using a specifically designed subroutine for making independent draws; the latter was used in all our examples.

2.3. *An illustration with a mixture CFTP.* Hobert, Robert and Titterington [HRT (1999)] considered the mixture model $pf_0(x) + (1-p)f_1(x)$ where $0 < p < 1$ is the estimand. Given i.i.d. observations $\{x_1, x_2, \ldots, x_n\}$ from this model, and a flat prior on $p$, the posterior for $p$ is proportional to $\prod_{i=1}^{n}[pf_0(x_i) + (1-p)f_1(x_i)]$. To deal with the continuous distribution of $p$, HRT use the natural discrete data augmentation that comes with the latent mixture indicator $\{z_1, z_2, \ldots, z_n\}$, where $z_i = 0$ if $x_i$ is from $f_0$ and $z_i = 1$ if $x_i$ is from $f_1$. Starting from $T = 1$, the following steps define the HRT algorithm:

1. Start a "bottom chain" at $p_{-T}^{(1)} = 0$ and a "top chain" at $p_{-T}^{(2)} = 1$.
2. At each $-T \leq t \leq -1$ and for $j = 1, 2$: generate $n$ i.i.d. uniform r.v.'s $\{u_{t1}, u_{t2}, \ldots, u_{tn}\} \equiv u_t$, and then set $z_{ti}^{(j)} = 0$ if $u_{ti} \leq \frac{p_t^{(j)} f_0(x_i)}{p_t^{(j)} f_0(x_i) + (1 - p_t^{(j)}) f_1(x_i)}$ and $z_{ti}^{(j)} = 1$ otherwise.

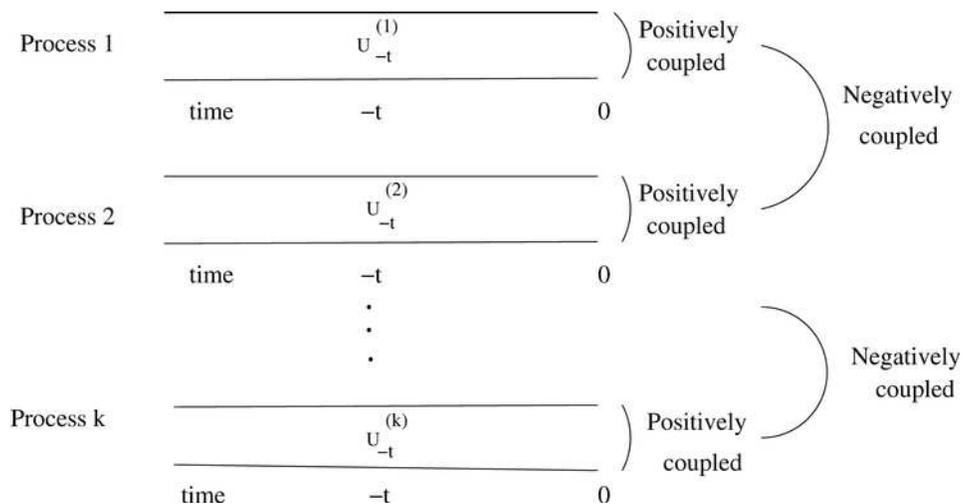

FIG. 1. *Parallel antithetic backward CFTP processes.*

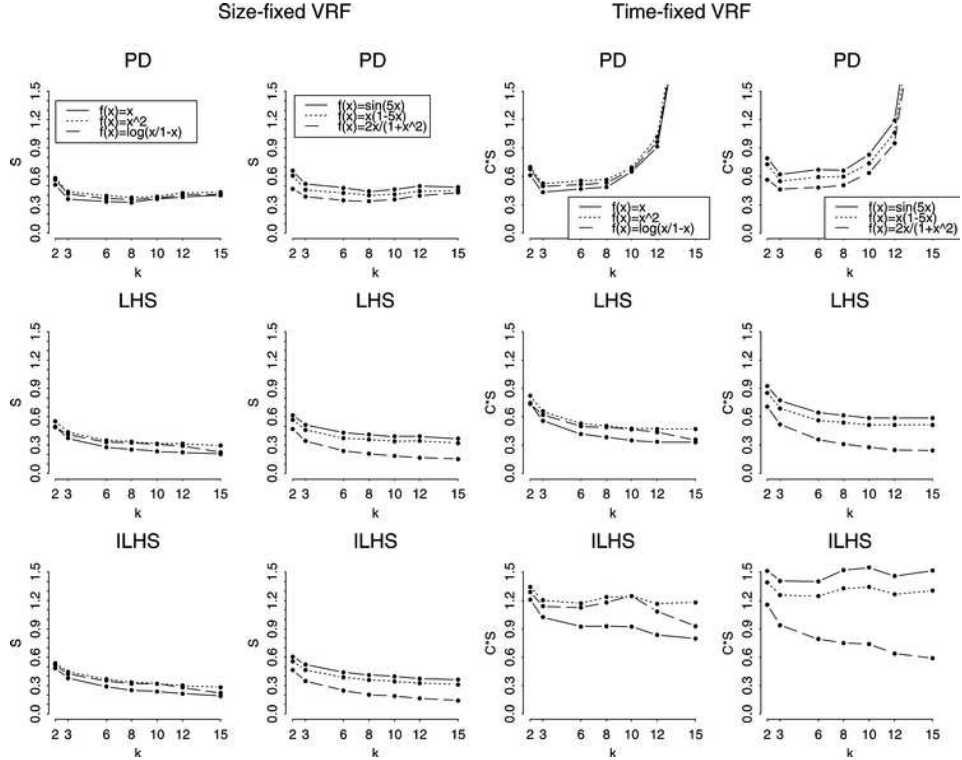

FIG. 2. *Normal mixture example. Size-fixed variance reduction factor* (left) *and time-fixed variance reduction factor* (right) *plotted against the number of parallel chains for different functions. Note the values of the time-fixed VRF for PD at $k = 15$ are too large for the plotting range.*

Let $m$ be the number of $z_{ti}^{(j)}$'s equal to 0 and let $w_t = \{w_{tr}, r = 1, \dots, n + 2\}$, where the $w_{tr}$'s are i.i.d. samples from an exponential distribution with mean 1. Compute $p_{t+1}^{(j)} = \psi(p_t^{(j)}, u_t, w_t) = \sum_{r=1}^{m+1} w_{tr} / \sum_{r=1}^{n+2} w_{tr}$.

3. If $p_0^{(1)} = p_0^{(2)} = p_0$, then $p_0$ is our sample. If not, set $T_{old} = T$ and go to step 2 with $T = 2T_{old}$ while keeping $u_{ti}, w_{tr}$ for all $-T_{old} \le t \le 0$, $1 \le i \le n$, $1 \le r \le n + 2$ unchanged.

Although $\psi(p, u, w)$ of step 2 is not monotone (in the same direction) in all its arguments, the variance reductions obtained with antithetic variates are still significant in all the simulation examples we have looked at. Figure 2 gives an illustration based on $0.33 \cdot \mathcal{N}(3.2, 3.2) + (1 - 0.33) \cdot \mathcal{N}(1.4, 1.4)$. The first two columns in Figure 2 plot the size-fixed VRF $S_k^{(f)}$ against $k$. The plots on each row correspond to a different generating method: the permuted displacement (PD) method, Latin hypercube sampling (LHS) and iterative Latin hypercube sampling (ILHS); see Section 4. The simulation size here is 7500 for each $k$, and the Monte Carlo variance for estimating $S_k^{(f)}$ is on the order of $10^{-6}$. We see that $S_k^{(f)}$ decreases from 0.5–0.6 with $k = 2$ to 0.2–0.3 when $k \ge 6$ (all numbers are from the ILHS method, which performs best). It appears that $S_k^{(f)}$ stabilizes after about $k = 10$. The second column shows the results for three nonmonotone estimand functions. The theoretical guarantee given in Section 3 does not apply to such functions, but nevertheless all $S_k^{(f)} < 1$. Even for $f(x) = \sin(5x)$, $S_k^{(f)}$ is around 0.4, implying a 60% reduction in variance when $k \ge 6$.

The last two columns in Figure 2 plot the time-fixed VRF, $T_k^{(f)}$ of (2.5). From the first row, it is seen that using $k \ge 10$ becomes too costly for PD, at least in our implementation. In contrast, the monotone patterns for the second and third rows resemble that for $S_k^{(f)}$, albeit the reduction is less because the $C_k$ factor tends to be more than 1. Nevertheless, we still see that $T_k^{(f)}$ reduces from 0.8 with $k = 2$ to about 0.5 when $k \ge 6$, and the empirical finding that using $k \ge 10$ is not practically beneficial remains true. These conclusions are based on the best performing method, LHS (ILHS is obviously more costly because of the iteration). However, we emphasize that in our implementation we have made no attempt to optimize our code. It is thus important to separate the gain in statistical efficiency, as measured by $S_k^{(f)}$, which does not depend on the particular implementation, from the possible offset, as measured by $C_k$, which depends critically on the particular implementation and thus could be further reduced with a more refined code/implementation.

2.4. *Forward antithetic coupling and slice sampling.* The next application used for illustration is a forward slice sampler. A graphical scheme of the antithetic principle applied to forward chains is shown in Figure 3. In this situation the parallel chains are started from different points in the sample space. After the so-called *burn-in period*, the realizations of each path are used, typically as positively correlated samples from the stationary distribution of the chain. Like CFTP, we can update $k$ chains using $k$ antithetically coupled variates. However, unlike CFTP, the within-chain autocorrelation for the forward implementation makes the determination of the VRF a bit more complicated. Specifically, we need to generalize (2.4) to

$$
\begin{aligned}
S_k^{(f)} &= 1 + \frac{\sum_{j_1 \neq j_2}^{k} \mathrm{Cov}(\sum_{t=1}^{m} f(X_t^{(j_1)}), \sum_{t=1}^{m} f(X_t^{(j_2)}))}{\sum_{j=1}^{k} \mathrm{Var}(\sum_{t=1}^{m} f(X_t^{(j)}))} \\
&= 1 + \frac{\sum_{j_1 \neq j_2}^{k} \sum_{t_1, t_2}^{m} \beta_{t_1, t_2; j_1, j_2}^{(f)}}{\sum_{j=1}^{k} \sum_{t_1, t_2}^{m} \gamma_{t_1, t_2; j}^{(f)}},
\end{aligned}
$$

(2.6)

where $\gamma_{t_1, t_2; j}^{(f)}$ and $\beta_{t_1, t_2; j_1, j_2}^{(f)}$ are respectively *within-chain* and *between-chain* autocovariance, that is,

$$
\begin{aligned}
\gamma_{t_1, t_2; j}^{(f)} &= \mathrm{Cov}(f(X_{t_1}^{(j)}), f(X_{t_2}^{(j)})), \\
\beta_{t_1, t_2; j_1, j_2}^{(f)} &= \mathrm{Cov}(f(X_{t_1}^{(j_1)}), f(X_{t_2}^{(j_2)})).
\end{aligned}
$$

(2.7)

In many applications expression (2.6) can be greatly simplified because of the ergodicity of the antithetically coupled joint chain (which is *not* an automatic consequence of the ergodicity of the marginal chains; see Section 3.2).

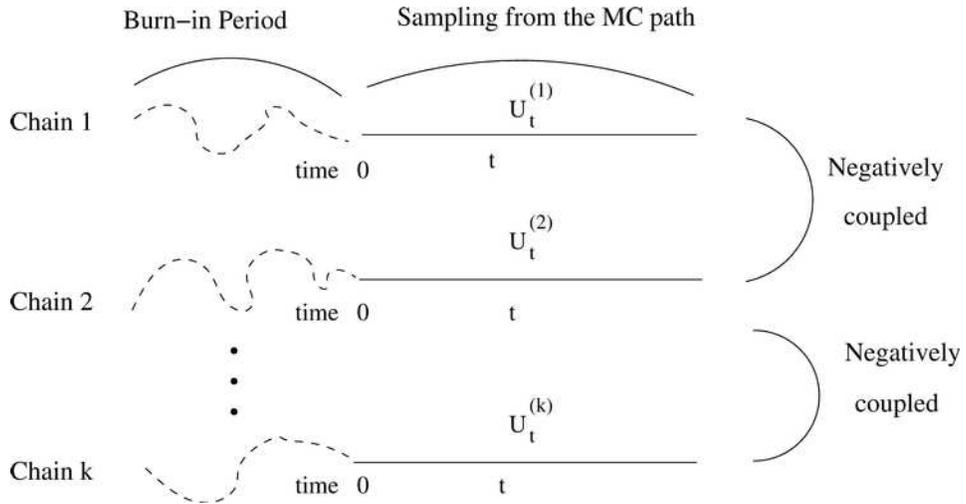

Fig. 3. *Parallel antithetic forward MCMC chains.*

Here we retain the full generality of (2.6) to show that as long as all the between-chain covariances are nonpositive, the use of antithetic coupling will always guarantee $S_k^{(f)} \leq 1$ regardless of the mixing properties of the individual chains (e.g., still in burning periods) or its within-covariance structure—that is, using antithetic coupling cannot hurt as long as $\beta_{t_1,t_2;j_1,j_2}^{(f)} \leq 0$ for all $(t_1, t_2, j_1, j_2)$; see Section 3 for conditions to guarantee this. The definition of time-fixed VRF remains the same as in (2.5) because $C_k$ is unaffected by the autocorrelation.

To illustrate possible gains, we adopt an example of Damien, Wakefield and Walker (1999), who use slice sampling [Neal (2003)] to draw from $\pi(x) \propto x^2 e^{-e^x} I_{\{x \geq 0\}}$. The resulting Gibbs sampler has the updating function,

$$(2.8) \qquad X_{t+1} = \psi(X_t, \xi_1, \xi_2) = \xi_1^{1/3} \log(e^{X_t} - \log(1 - \xi_2)),$$

where $\xi_1$ and $\xi_2$ are i.i.d. Uniform$(0, 1)$. Since $\psi(x, \xi_1, \xi_2)$ is nondecreasing, Theorem 1 of Section 3 guarantees that $S_k^{(f)} \leq 1$ for any monotone function $f$. Note (2.8) is an example where $\psi$ is *made* to be nondecreasing by using $1 - \xi_2$ instead of $\xi_2$ inside the second logarithm.

Figure 4 is the counterpart of Figure 2 for the forward case, with similar general patterns (the simulation size here is 5000). In particular, for monotone functions $S_k^{(f)}$ decreases from between 0.35–0.45 with $k = 2$ to 0.1–0.15 when $k \geq 6$. For highly nonmonotone $f(x) = \sin(5x)$, in contrast to $k = 2$, $S_k^{(f)} \leq 1$ when $k \geq 3$. For less variable nonmonotone functions $f(x) = 2x(1 + x^2)^{-1}$ and $f(x) = x(1 - 5x)$, $S_k^{(f)}$ goes below 0.5. Partial theoretical support for this phenomenon can be found in Owen (1997), whose result implies that, under LHS, $S_k^{(f)} \leq k/(k-1)$ for any square-integrable $f$. This suggests that there is more room for $S_2^{(f)}$ to exceed 1 than for $S_k^{(f)}$ with $k > 2$. However, in terms of the time-fixed VRF, the clear winner is the LHS method with monotone functions. For highly nonmonotone functions, our implementation becomes too costly, suggesting that if such estimand functions are of main interest, then it would be generally safer to just use the independent implementation of multiple chains, as recommended in Gelman and Rubin (1992).

2.5. *A real-data application with Bayesian probit regression.* We conclude our empirical investigation by presenting a multidimensional real-data application. The data are taken from van Dyk and Meng (2001) and consist of measurements on 55 patients, of whom 18 have been diagnosed with latent membranous lupus. Table 1 shows the data with two clinical covariates, IgA and IgG, that measure the levels of immunoglobulin of type A and of type G, respectively. Of interest is the prediction of disease occurrence using the two covariates IgG3 − IgG4 and IgA. As in van Dyk and Meng (2001), we consider a probit regression model. For each patient $i$, we model the disease indicator variables as independent $Y_i \sim \text{Bernoulli}(\Phi(x_i^T \beta))$, where $\Phi(\cdot)$ is the CDF of $N(0, 1)$, $x_i$ is the vector of the covariates and $\beta$ is a $3 \times 1$ vector of parameters. We want to sample from the posterior distribution corresponding to the flat prior for $\beta$.

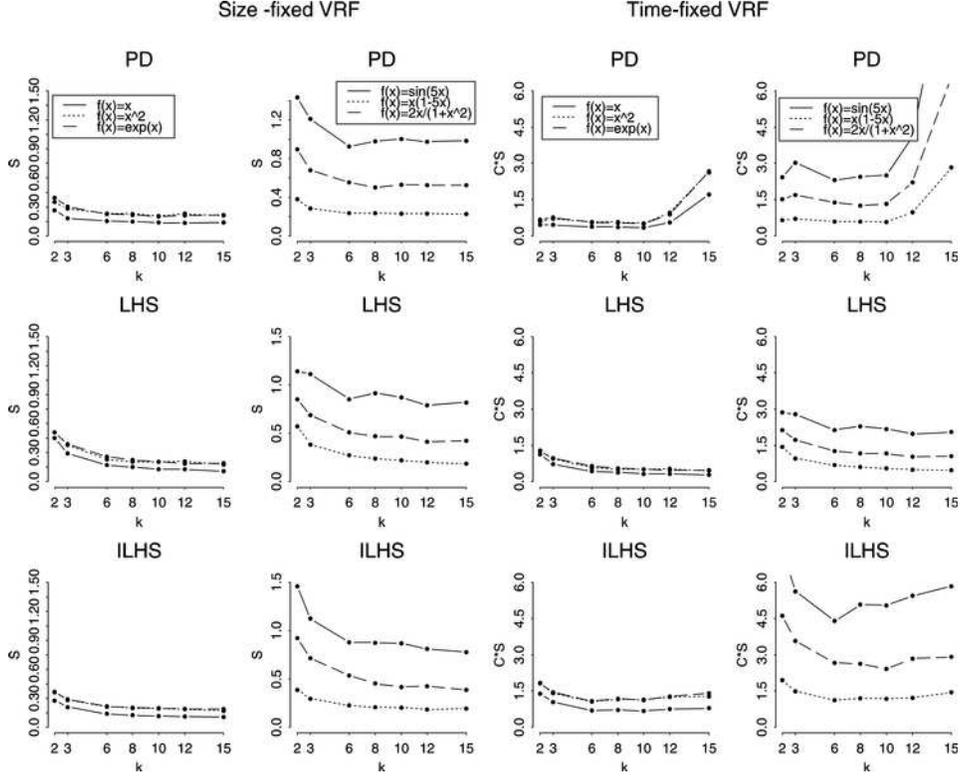

Fig. 4. *Slice sampling example. Size-fixed variance reduction factor* (left) *and time-fixed variance reduction factor* (right) *plotted against the number of parallel chains for different functions. Note the change of scale in the vertical axis from the sized-fixed VRF to time-fixed VRF.*

To illustrate the impact of antithetic coupling on mixing, we adopt the standard Gibbs sampler with the latent variables $\psi_i \sim N(x_i^T \beta, 1)$, of which only the sign $Y_i$ is observed, as the augmented data [e.g., Albert and Chib (1993)]. Let $X$ be the $n \times p$ matrix whose $i$th row is $x_i$ and $\psi = (\psi_1, \ldots, \psi_n)$; for our example $n = 55$ and $p = 3$. The resulting Gibbs sampler alternates between sampling from (i) $\beta | \psi \sim N(\tilde{\beta}, (X^T X)^{-1})$ with $\tilde{\beta} = (X^T X)^{-1} X^T \psi$ and from (ii) $\psi_i | \beta, Y_i \sim TN(x_i^T \beta, 1, Y_i)$, where $TN(\mu, \sigma^2, Y)$ is $N(\mu, \sigma^2)$ truncated to be positive if $Y > 0$ and negative if $Y < 0$. Our parallel chains are then coupled antithetically at each of these two updating steps.

TABLE 1
*The number of latent membranous lupus nephritis cases, the numerator, and the total number of cases, the denominator, for each combination of the values of the two covariates*

| IgG3 − IgG4 | IgA | | | | |
|:---:|:---:|:---:|:---:|:---:|:---:|
| | **0** | **0.5** | **1** | **1.5** | **2** |
| −3.0 | 0/1 | — | — | — | — |
| −2.5 | 0/3 | — | — | — | — |
| −2.0 | 0/7 | — | — | — | 0/1 |
| −1.5 | 0/6 | 0/1 | — | — | — |
| −1.0 | 0/6 | 0/1 | 0/1 | — | 0/1 |
| −0.5 | 0/4 | — | — | 1/1 | — |
| 0 | 0/3 | — | 0/1 | 1/1 | — |
| 0.5 | 3/4 | — | 1/1 | 1/1 | 1/1 |
| 1.0 | 1/1 | — | 1/1 | 1/1 | 4/4 |
| 1.5 | 1/1 | — | — | 2/2 | — |

In the first step we want to generate exchangeable $\{\beta^{(1)}, \ldots, \beta^{(k)}\}$ with $\beta^{(j)} \in \mathbf{R}^p$ for all $1 \le j \le k$, such that marginally $\beta^{(j)} \sim N(\tilde{\beta}, \Sigma = (X^T X)^{-1})$ and $\mathrm{Corr}(\beta_i^{(j)}, \beta_i^{(j')})$ is minimum possible, for all $1 \le i \le p$ and all $1 \le j \ne j' \le k$. This can be realized by generating $p$ i.i.d. multivariate normal $(Z_i^{(1)}, \ldots, Z_i^{(k)})^\top$ such that $Z_i^{(j)} \sim N(0,1)$ and $\mathrm{Corr}(Z_i^{(j)}, Z_i^{(j')}) = -1/(k-1)$. We then let $\beta^{(j)} = \Sigma^{1/2} Z^{(j)} + \tilde{\beta}$, where $Z^{(j)} = (Z_1^{(j)}, \ldots, Z_p^{(j)})^\top$ and $\Sigma^{1/2}$ is the Choleski decomposition of $\Sigma$.

In the second step, we use the inverse CDF method suggested by Gelfand, Smith and Lee ([1992](#)) to sample from $TN(\mu, \sigma^2, 1)$. Namely, we simulate $U \sim \mathrm{Uniform}(0,1)$ and then take

$$(2.9) \qquad Z = \mu + \sigma \Phi^{-1}[\Phi(-\mu/\sigma) + U(1 - \Phi(-\mu/\sigma))].$$

The methods of Lew ([1981](#)) and Bailey ([1981](#)) were used to approximate $\Phi(x)$ and $\Phi^{-1}(x)$ needed for ([2.9](#)). The antithetic coupling is then realized via $n$ independent vectors of NA uniform random variables $\{U_i^{(1)}, \ldots, U_i^{(k)}; 1 \le i \le n\}$ so that for the $j$th chain $\psi_i^{(j)}$ is generated by using $U = U_i^{(j)}$ in ([2.9](#)). We used only two iterations of ILHS, to increase the speed of the antithetic sampler. The computational overhead is small even in the first step as all the matrix inversions required there were performed once outside the inner loop of the stochastic algorithm.

The improvement brought about by the antithetic coupling is twofold because it can reduce both the variance and the bias of the original sampler. Van Dyk and Meng ([2001](#)) demonstrate the slow mixing of the standard algorithm and propose a much more reliable marginal data augmentation algorithm. Our antithetic implementation does not, and cannot, remedy all the problems with the original sampler, but it can provide noticeable improvement without requiring a new algorithm or a substantial increase in the computational load. This is seen in Figure 5, which contains scatter plots of $(\beta_0, \beta_1)$ using draws from $k$ independent chains of the original sampler (left column) and from $k$ antithetic chains (right column). The true contour plots, obtained via numerical methods, are from van Dyk and Meng (2001). Each scatter plot is based on 9000 sample points, divided equally among the $k$ chains used. All the starting points were set to the same MLE; we did this on purpose to see to what extent the antithetic coupling can "spread out" despite the fact that all chains are started at the same point. Although the antithetic chains still have the serious problem of missing a good part of the "right tail," a problem that was avoided by van Dyk and Meng's (2001) algorithm, in all rows the scatter plot in the right column extends to the "right tail" no less, and sometime substantially more, than the scatter plot in the left column.

To represent quantitatively the benefit of antithetic coupling, Figure 6 plots, against $k$, the relative bias, relative standard error (SD) and relative root mean square error (RMSE) from antithetic chains, in estimating six posterior quantities, with the results from independent chains as the baseline (e.g., a relative bias 0.4 means the antithetic implementation reduces the bias, in magnitude, by $1 - 0.4 = 60\%$). The six quantities we choose are, in the order of the rows in Figure 6, $\mathrm{E}[\beta_0|\mathcal{D}]$, $\mathrm{E}[\beta_1|\mathcal{D}]$, $\mathrm{Var}[\beta_0|\mathcal{D}]$, $\mathrm{Var}[\beta_1|\mathcal{D}]$, $\mathrm{E}[-\beta_0/\beta_1|\mathcal{D}]$ and $\mathrm{E}[Q|\mathcal{D}]$, where $Q = [\Phi(\beta_0 - 0.5\beta_2 + 1.5\beta_2)]/[1 - \Phi(\beta_0 - 0.5\beta_2 + 1.5\beta_2)]$ and $\mathcal{D}$ denotes the observed data $(X, Y)$. Whereas the first four are usually required in a Bayesian analysis, the last two are specific to the example at hand: $-\beta_0/\beta_1$ is the so-called LD50 level (i.e., with 50% chance the corresponding dosage becomes lethal) for covariate $x_1$ when $x_2 = 0$, and $Q$ represents the odds of having the disease when $x_1 = -0.5$ and $x_2 = 1.5$. In particular, four out of these six estimand functions are nonmonotone.

The simulation sizes in Figure 6 are the same as in Figure 5, and the reduction factors are then computed using 1000 replicates of the simulation process. No burn-in period was discarded because we are interested in comparing the mixing properties of the two implementations. Because the total sample size is fixed at 9000 for each plot, using a larger $k$ means a smaller within-chain sample size $m$, and hence possibly more strong influence of the starting point. To investigate the effect of the starting point, we use three starting strategies: 1. MLE (columns 1–3 in Figure 6). 2. At two standard deviations (SD) from the MLEs, in this case, $\beta_0 + 2SD$ and $\beta_i - 2SD$, $i = 1, 2$ (columns 4–6). 3. From extreme points (more than three SD away from the MLE) in the support of the distribution (columns 7–9). In general it is seen that the RMSE is reduced due to a decrease in variance, bias or both. Best results correspond to $3 \le k \le 6$, reflecting the aforementioned trade-off between $k$ and $m$. As expected, the reduction factors vary with the starting points, as well as the estimand of interest. Nevertheless, with $k \ge 3$ the antithetic coupling in general does not inflate the RMSE and may result in reduction factors as low as 0.1, that is, up to 90% of saving.

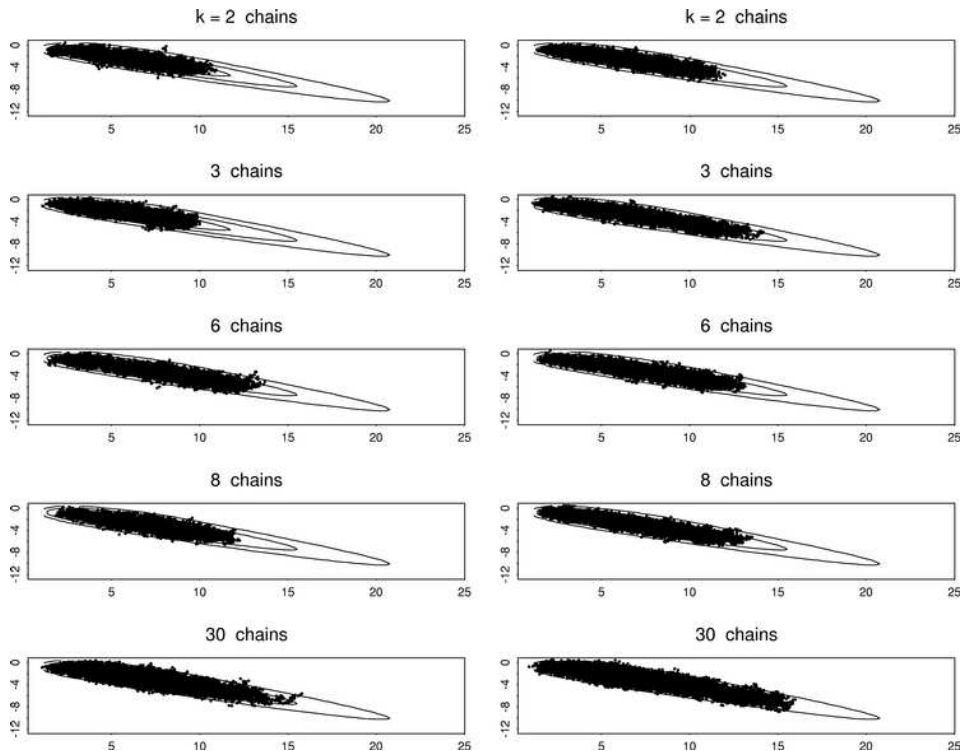

FIG. 5. *Bayesian probit regression example. Scatter plots of* $(\beta_0, \beta_1)$ *generated by using* $k$ *independent chains* (left) *and* $k$ *antithetic chains* (right) *for different values of* $k$; *the contour lines are from the targeted bivariate posterior distribution. Each plot contains* 9000 *draws.*

## 3. A theoretical foundation.

3.1. *Negative association and dependence.* In general, a qualitative measure of negative dependence should adequately reflect the following intuitive behavior among a set of variables: if one subset of the variables is "high," then a disjoint subset of the variables is "low." Different ways to define such negative dependence have received a great deal of attention in the last twenty years or so. Due to the success of the *positive association* concept of Esary, Proschan and Walkup (**1967**), the main challenge has been to build the negative association concepts as "duals" of the positive ones, but so far there has not been a universally acceptable construction [e.g., Pemantle (2000)]. Specifically, the set of random variables $\{X_1, \ldots, X_n\}$ is said to be *positively associated* (PA) if for any nondecreasing (or nonincreasing—we will not state both hereafter) functions $f_1, f_2$, $\mathrm{Cov}(f_1(X_1, X_2, \ldots, X_n), f_2(X_1, X_2, \ldots, X_n)) \geq 0$. Here a function $f: \mathbf{R}^n \to \mathbf{R}$ is called nondecreasing if it is nondecreasing in all of its arguments.

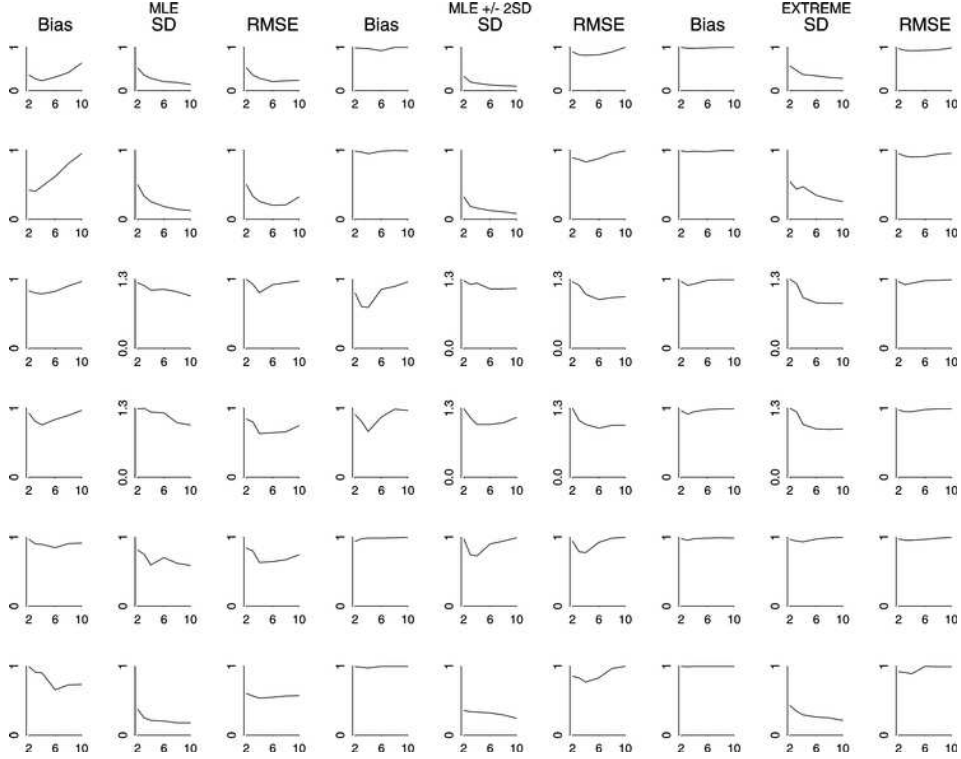

Fig. 6. *Relative Monte Carlo bias, standard error (SD) and root mean squared error (RMSE), from the antithetic chains, relative to that of independent chains, as functions of $k$, the number of parallel chains. Same simulation configurations as in Figure 5, and with 1000 replication for each $k$. Each row represents a different estimand function; from top to bottom:* $\mathrm{E}[\beta_0|\mathcal{D}]$; $\mathrm{E}[\beta_1|\mathcal{D}]$; $\mathrm{Var}(\beta_0|\mathcal{D})$; $\mathrm{Var}(\beta_1|\mathcal{D})$; $\mathrm{E}[-\beta_0/\beta_1|\mathcal{D}]$; $\mathrm{E}[Q|\mathcal{D}]$. *Columns 1–3 use MLE as the starting points, 4–6 use MLE $\pm 2$ SD, and 7–9 use some extreme points situated more than three SD away from the MLE.*

The closest negative counterpart of positive association we can find is the *negative association* concept introduced by Joag-Dev and Proschan (1983).

DEFINITION 1. The random variables $\{X_i\}_{1 \leq i \leq n}$, where each $X_i$ can be of arbitrary dimension, are said to be negatively associated (NA) if for every pair of disjoint finite subsets $A_1, A_2$ of $\{1, 2, \ldots, n\}$ and for any nondecreasing

functions $f_1$, $f_2$

$$(3.1) \qquad \mathrm{Cov}(f_1(X_i, i \in A_1), f_2(X_j, j \in A_2)) \leq 0$$

whenever the above covariance function is well defined.

The following equivalence result is important for some of our subsequent theoretical investigation.

PROPOSITION 1. *The random variables* $\{X_i\}_{1 \leq i \leq n}$ *are NA if and only if* (3.1) *holds for every pair of nonnegative, bounded and nondecreasing functions* $f_1$ *and* $f_2$.

PROOF. The "only if" part is obvious. To prove the "if" part, let $f_m(X_i, i \in A_m)$, $m = 1, 2$, be the two functions in (3.1). For any positive integer $l$, let $f_m^{(l)}$ be the truncation of $f_m$ to $[-l, l]$, that is, $f_m^{(l)} = f_m$ when $|f_m| \leq l$, and $f_m^{(l)} = \pm l$ depending on $f_m > l$ or $f_m < -l$. Clearly, $|f_m^{(l)}| \leq |f_m|$, and thus by the dominated convergence theorem, $\lim_{l \to \infty} \mathrm{Cov}(f_1^{(l)}(X_i, i \in A_1), f_2^{(l)}(X_j, j \in A_2)) = \mathrm{Cov}(f_1(X_i, i \in A_1), f_2(X_j, j \in A_2))$. This allows us to conclude (3.1) because

$$\mathrm{Cov}(f_1^{(l)}(X_i, i \in A_1), f_2^{(l)}(X_j, j \in A_2))$$
$$= \mathrm{Cov}(f_1^{(l)}(X_i, i \in A_1) + l, f_2^{(l)}(X_j, j \in A_2) + l) \leq 0,$$

where the last inequality holds because $f_m^{(l)} + l$ is a nonnegative, bounded and nondecreasing function for $m = 1, 2$. □

The notion of NA is most useful for our purposes primarily because, like PA, it is closed under the independent union operation as well as monotone transformations, as proved in Joag-Dev and Proschan (1983). Specifically, we have:

PROPOSITION 2. *If* $\{X_1, \ldots, X_{n_1}\}$ *and* $\{Y_1, \ldots, Y_{n_2}\}$ *are two independent sets of NA (PA) random variables, then their union,* $\{X_1, \ldots, X_{n_1}, Y_1, \ldots, Y_{n_2}\}$, *is a set of NA (PA) random variables.*

PROPOSITION 3. *If* $\{X_i\}_{1 \leq i \leq n}$ *is a sequence of NA (PA) random variables and* $(\psi_i)_{1 \leq i \leq n}$ *are all nondecreasing functions, then* $(\psi_i(X_i))_{1 \leq i \leq n}$ *is a sequence of NA (PA) random variables.*

These two results enable us to prove the following fundamental result on antithetic coupling of homogeneous and nonhomogeneous Markov chains.

Theorem 1. *Suppose we run a $k$-process antithetically coupled chain, $\mathcal{X}_t^k = \{X_t^{(1)}, \ldots, X_t^{(k)}\}$, forward for $T$ iterations, where $X_t^{(j)} = \psi_t(X_{t-1}^{(j)}, U_t^{(j)})$, $j = 1, \ldots, k$, and $\mathcal{U}_t^k = \{U_t^{(1)}, \ldots, U_t^{(k)}\}$, $t = 1, \ldots, T$, are $T$ independent sets of NA variables. Assuming $\psi_t(X, U)$ is nondecreasing for all $t \leq T$, we have the following results.*

(i) *The $k$-tuple $\{X_{t_1}^{(1)}, \ldots, X_{t_k}^{(k)}\}$ is a collection of $k$ NA variables for any $\{t_1, \ldots, t_k\} \in \{0, 1, \ldots, T\}^k$ if and only if it is so for $t_1 = \cdots = t_k = 0$.*

(ii) *Assuming $\mathcal{X}_0^k$ is a set of NA variables, then the variance reduction factor, $S_k^{(f)}$, defined by (2.6), is at most 1 for any monotone function $f$, and it is strictly less than 1 if and only if at least one of the between-chain covariances, $\beta_{t_1, t_2; j_1, j_2}^{(f)}$ of (2.7), is strictly negative.*

Proof. For (i) the necessity holds by definition. For the sufficiency, because $\{X_0^{(1)}, \ldots, X_0^{(k)}\}$ and $\{U_t^{(1)}, \ldots, U_t^{(k)}\}$, $t = 1, \ldots, T$, are $T + 1$ sets of independent NA variables, by Proposition 2 their union is also NA. Consequently, because $X_{t_j}^{(j)}$ is a nondecreasing function of $Z^{(j)} \equiv \{X_0^{(j)}, U_t^{(j)}, t = 1, \ldots, T\}$ only, and $Z^{(i)}$ and $Z^{(j)}$ do not share any common variable for $i \neq j$, by Proposition 3 $\{X_{t_1}^{(1)}, \ldots, X_{t_k}^{(k)}\}$ are NA for any $t_j \leq T$, $j = 1, \ldots, k$.

For (ii), by (i) all the between-chain autocovariances $\beta_{t_1, t_2; j_1, j_2}^{(f)} \leq 0$, which implies $S_k^{(f)} \leq 1$ because the denominator in (2.6), being a variance, is always positive. (In fact, by using the PA part of Propositions 2 and 3, we can conclude all $\gamma_{t_1, t_2; j}^{(f)} \geq 0$.) Furthermore, because all $\beta_{t_1, t_2; j_1, j_2}^{(f)} \leq 0$, by (2.6), $S_k^{(f)} = 1$ if and only if all $\beta_{t_1, t_2; j_1, j_2}^{(f)}$'s are zero. □

Since for a Gibbs sampler with attractive stationary density the updating function can be expressed as a monotone function [e.g., Propp and Wilson (1996) and Häggström and Nelander (1998)], Theorem 1 covers Frigessi, Gåsemyr and Rue's (2000) Theorem 1 with $k = 2$. For practical purposes, the requirement that $\mathcal{X}_0^k = \{X_0^{(1)}, \ldots, X_0^{(k)}\}$ are NA is immaterial, because being fixed (even with the choice that $X_0^{(1)} = \cdots = X_0^{(k)}$) or more generally being independent is a trivial case of being NA. It is also evident that Theorem 1 does not require $X_0^{(j)}$ to be from the stationary distribution.

Theorem 1, however, does require that the random variables in the $k$-tuple $\mathcal{U}_t^k = \{U_t^{(1)}, \ldots, U_t^{(k)}\}$ are NA. Typically, for a particular distribution there are many ways of achieving this, and some general recipes are described in Section 4. A useful construction in practice is to use the fact that monotone functions of distinctive subsets of NA variables are NA, a fact that allows us to build upon known NA variables, such as permutation distribu-tions, multinomial, multivariate hypergeometric, Dirichlet and multivariate normals with nonpositive correlations [see Joag-Dev and Proschan (1983)].

We remark that part (i) of Theorem 1 is actually stronger than we need, because from (ii), in order to assure $S_k^{(f)} \leq 1$, we only need to ensure pairwise negative covariance: $\mathrm{Cov}(f(X_{t_1}^{(j_1)}), f(X_{t_2}^{(j_2)})) \leq 0$. This weaker version is sometimes easier to establish, and thus we define the following notion.

DEFINITION 2. The random variables $X_1, X_2, \ldots, X_n$ are said to be pairwise negatively associated (PNA) if $\{X_i, X_j\}$ are NA for any $i \neq j \in \{1, \ldots, n\}$.

Clearly, Propositions 2 and 3 still hold for PNA. In the next section, we will see that only requiring our outcome to be PNA can avoid technical complications that are not critical to MCMC applications in practice. In particular, the PNA property is easier to verify because it is equivalent to the *negative quadrant dependency* (NQD) notion for *a pair* of random variables, as defined by Lehmann (1966). For more than two variables, NQD is a consequence of NA, but not vice versa, as formalized by the following result, also due to Joag-Dev and Proschan (1983).

PROPOSITION 4. *Let $X_1, X_2, \ldots, X_n$ be a set of NA random variables. Then they are also negatively lower orthant dependent (NLOD), that is,*

$$(3.2) \quad P(X_1 \leq x_1, \ldots, X_n \leq x_n) \leq \prod_{i=1}^n P(X_i \leq x_i) \qquad \text{for all } x_1, \ldots, x_n,$$

*as well as negatively upper orthant dependent (NUOD), namely,*

$$(3.3) \quad P(X_1 > x_1, \ldots, X_n > x_n) \leq \prod_{i=1}^n P(X_i > x_i) \qquad \text{for all } x_1, \ldots, x_n.$$

3.2. *Negative association in limit and joint convergence.* Although Theorem 1 provides a theoretical foundation for antithetic coupling with forward MCMC, it does not cover backward MCMC such as CFTP. This is because the $T$ in Theorem 1 is a nonrandom constant, whereas for CFTP it is a random variable and, more importantly, it is not independent of the corresponding draw from the CFTP. We thus need to extend Theorem 1 to the case with $T = +\infty$, namely, we need to prove that the limiting $k$-tuple $\{X_\infty^{(1)}, \ldots, X_\infty^{(k)}\}$ is still NA, or at least PNA. This would entail that the joint draw from a CFTP, $\{X_0^{(1)}, \ldots, X_0^{(k)}\}$, is NA/PNA because we can identify its probability structure with that of $\{X_\infty^{(1)}, \ldots, X_\infty^{(k)}\}$ from a forward process.

The extension to $T = \infty$ is not straightforward because the ergodicity of the marginal chain does not guarantee that of the joint chain, so that a limiting argument such as $\text{Cov}(f(X_\infty^{(i)}), f(X_\infty^{(j)})) = \lim_{t \to \infty} \text{Cov}(f(X_t^{(i)}), f(X_t^{(j)})) \le 0$ needs qualification. The problem has been discussed also by Arjas and Gasbarra (1996) and Frigessi, Gåsemyr and Rue (2000) who have shown that, if the joint chain is $\phi$-irreducible and the closure of the support of $\phi$ contains an open set, then the joint chain is positive recurrent on the closure of the support of $\phi$. This can be used to prove that the above convergence of covariances holds. However, the assumption that the support of the irreducibility measure $\phi$ has nonempty interior is violated even by simple examples such as the following Markov chain: take $X_0 \in S_1$, the unit circle; for $t \ge 1$, draw $\theta_t$ uniform on $[0, \pi)$ and construct the line that goes through $X_{t-1}$ and has slope $\tan(\theta_t)$. The intersection of this line with the unit circle is $X_t$. The left panel of Figure 7 illustrates this construction with $t = 1$. The right panel illustrates a pair of antithetically coupled chains $(X_t, Y_t)$ (for $t \le 2$) where the $Y_t$ is constructed the same way as $X_t$ except each time $X_t$ is updated using $\theta_t$, $Y_t$ is updated using $\pi - \theta_t$. Algebraically, we have

$$(3.4) \quad \begin{aligned} \theta_t^X &= (-\theta_{t-1}^X + \pi + 2\theta_t) \bmod (2\pi), \\ \theta_t^Y &= (-\theta_{t-1}^Y + \pi - 2\theta_t) \bmod (2\pi), \end{aligned}$$

where $\theta_t^X$ and $\theta_t^Y$ are the polar-coordinate representations of $X_t$ and $Y_t$, respectively. Marginally, $\theta_t^X$ (and hence $\theta_t^Y$) are i.i.d. uniform variables on $[0, 2\pi)$, and hence $\{X_t\}_t$ is trivially uniformly ergodic. However, the joint state space of $(\theta_t^X, \theta_t^Y)$, $S_2 = [0, 2\pi) \times [0, 2\pi)$, consists of *uncountably* many

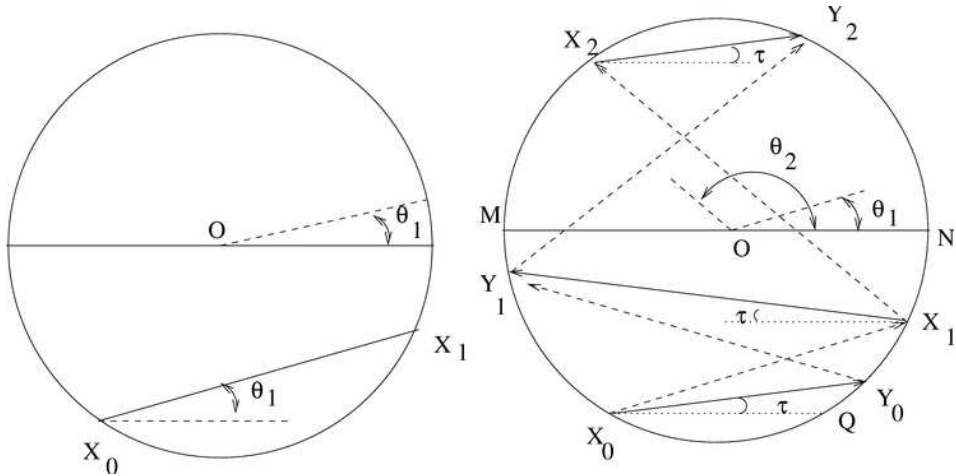

Fig. 7. *The random walk on the circle example. Illustration of construction of $X_t$ (left) and of the antithetic coupling $(X_t, Y_t)$ (right).*

absorbing subclasses, defined by

$$(3.5) \quad S(\tau) = \{(\theta^X, \theta^Y) \in S_2 : \theta^X + \theta^Y = (2k+1)\pi \pm 2\tau, \text{ for } k = 0, 1\},$$

where $\tau \in [0, \pi/2]$ is the acute angle between the line $\overline{X_0 Y_0}$ and the horizontal diameter (see Figure 7). Here $S(\tau)$ is an absorbing set because the joint updating rule leaves the acute angle between the line $\overline{X_t Y_t}$ and the horizontal diameter, denoted by $\tau_t$, unchanged. Let $\rho_t(\tau) = \text{Corr}(X_t^{(1)}, Y_t^{(1)})$, where $Z^{(1)}$ denotes the $X$-axis coordinate of $Z$ [i.e., $Z^{(1)} = \cos(\theta^Z)$]. Then, as proved in Appendix A.1, $\rho_t(\tau) = -\cos(2\tau)$, increasing from $-1$ to $1$ as $\tau$ increases from $0$ to $\pi/2$. In addition, for each of the ergodic classes, the joint chain has an irreducibility measure whose support has an empty interior.

This example shows that the uniform ergodicity of the marginal chain is not sufficient to guarantee the ergodicity of the antithetically coupled chain or to directly justify the extension of the NA properties from $\{X_t^{(1)}, \ldots, X_t^{(k)}\}$ to its "limit" $\{X_\infty^{(1)}, \ldots, X_\infty^{(k)}\}$. Here the "limit" is in quotes because while the marginal distribution of $X_\infty^{(j)}$ does not depend on the starting point, the joint distribution of $\{X_\infty^{(1)}, \ldots, X_\infty^{(k)}\}$ may depend on it, as is the case in the previous example. Fortunately, the time-backward dual sequence approach discussed in Section 2.1 provides a theoretical tool to bypass the issue of the joint chain's ergodicity, and thereby to establish the following counterpart of Theorem 1 for justifying backward antithetic coupling for all monotone CFTP algorithms.

THEOREM 2. *Let $X_t = \psi(X_{t-1}, U_t)$ be a uniformly ergodic Markov chain on state space $\Gamma$ with $\pi$ being its invariant distribution and with $\psi$ nondecreasing. Let $\mathcal{X}_t^k = \{X_t^{(1)}, \ldots, X_t^{(k)}\}$ be the corresponding antithetically coupled joint chain as in Theorem 1. Then:*

(i) *For any given starting point $\mathcal{X}_0^k$, $\mathcal{X}_t^k = \{X_t^{(1)}, \ldots, X_t^{(k)}\}$ converges in distribution to some $\mathcal{X}_\infty^k = \{X_\infty^{(1)}, \ldots, X_\infty^{(k)}\}$, whose joint distribution may depend on $\mathcal{X}_0^k$.*

(ii) *For any nondecreasing real functions $f_m$, $m = 1, 2$, on $\Gamma$ such that $\text{Var}_\pi(f_m(X)) < \infty$, $m = 1, 2$, and such that their discontinuity points are contained in a $\pi$-null set, we have*

$$(3.6) \quad \text{Cov}(f_1(X_\infty^{(j_1)}), f_2(X_\infty^{(j_2)}) | \mathcal{X}_0^k) \leq 0 \qquad \text{for any } j_1 \neq j_2.$$

PROOF. For each marginal chain $X_t^{(j)}$, let $\tilde{X}_t^{(j)}$ be its dual sequence as in (2.2). Then, by the equivalence of the implementability of CFTP and uniform ergodicity [Foss and Tweedie (1998)], we know that there exists an almost sure finite stopping time $S$ on the infinite-product space $\Gamma_U^\infty = \prod_{t \geq 1} \mathcal{U}_t^k$ such that $\tilde{X}_S^{(j)} \sim \pi$ and $\tilde{X}_t^{(j)} = \tilde{X}_S^{(j)}$ for all $t \geq S$ and all $j$. Therefore, by re-expressing all relevant variables on their common sampling space $\Gamma_U^\infty$, we can define $\mathcal{X}_\infty^k(\omega) = \{\tilde{X}_{S(\omega)}^{(1)}(\omega), \ldots, \tilde{X}_{S(\omega)}^{(k)}(\omega)\}$, which is a well-defined $k$-component joint random variable on $\Gamma_U^\infty$ because $P(S(\omega) < \infty) = 1$. Furthermore, because $\tilde{\mathcal{A}}_t^k = \{\tilde{X}_t^{(1)}, \ldots, \tilde{X}_t^{(k)}\}$ converges to $\mathcal{X}_\infty^k$ with probability 1 by the construction, and because $\mathcal{A}_t^k$ and $\tilde{\mathcal{A}}_t^k$ have the same distribution for any $t$, $\mathcal{A}_t^k$ must converge *in distribution* to $\mathcal{X}_\infty^k$, and hence (i).

For (ii), we only need to prove (3.6) when both $f$'s are bounded, following the same argument as in the proof of Proposition 1. Since $f_m, m = 1, 2$, are almost surely continuous with respect to the distribution of $\{X_\infty^{(i)}, X_\infty^{(j)}\}$ because its margins are $\pi$, we can conclude by part (i) that $\{f_1(X_t^{(i)}), f_2(X_t^{(j)})\}$ converges in distribution to $\{f_1(X_\infty^{(i)}), f_2(X_\infty^{(j)})\}$. This implies

$$(3.7) \quad \mathrm{Cov}(f_1(X_\infty^{(i)}), f_2(X_\infty^{(j)})|\mathcal{A}_0^k) = \lim_{t\to\infty} \mathrm{Cov}(f_1(X_t^{(i)}), f_2(X_t^{(j)})|\mathcal{A}_0^k) \leq 0,$$

where the limiting argument holds because both $f$'s are bounded. $\square$

Part (ii) of Theorem 2 shows that $\{\tilde{X}_\infty^{(1)}, \ldots, \tilde{X}_\infty^{(k)}\}$ are PNA when the monotone functions are almost surely continuous with respect to the underlying dominating measure, which is the case if the latter is the Lebesgue measure, as is common in Bayesian computation. We emphasize that the condition that $\psi$ is monotone is not needed for part (i), but it is important for part (ii). This is seen in the unit circle example, where (i) holds with uncountably many limiting distributions for the joint chain, depending on the initial $\tau$. However, (ii) does not hold because the $\psi$ function there is obviously not monotone. The fact that $\rho(\tau) < 0$ when $\tau < \pi/4$ in that example also illustrates that the monotonicity is a sufficient but not necessary condition for preserving negative correlation.

3.3. *Extreme antithesis.* The concept of NA provides a qualitative description of negative dependence. Quantitatively, it is desirable to generate $\{X_1, \ldots, X_k\}$ (corresponding to $\{X^{(1)}, \ldots, X^{(k)}\}$ in previous sections) such that $\mathrm{Corr}(f(X_i), f(X_j))$ is as negative as possible. Formally, we define the following notion.

DEFINITION 3. A set of variables $\{X_1, X_2, \ldots, X_k\}$ is said to achieve extreme antithesis (EA) with respect to a (marginal) distribution $F$ if they are exchangeable and

$$\mathrm{Corr}(X_i, X_j) = \min\{\mathrm{Corr}(Y_i, Y_j) : Y_1, \ldots, Y_k \text{ exchangeable}, Y_i \sim F \ \forall i\}.$$

For $k = 2$, the single strategy of using quantile coupling via $X_1 = F^{-1}(U)$ and $X_2 = F^{-1}(1 - U)$ achieves EA for any $F$, as discussed in Section 1. For $k > 2$, the matter is much more complicated, even just for $F :=$ Uniform$(0, 1)$ and $F := N(0, 1)$. Indeed, it is not hard to establish the following negative result [see Craiu and Meng (2002) for a proof], where $\Phi$ is the CDF of $N(0, 1)$.

PROPOSITION 5. *It is impossible to find a joint distribution* $F^{(3)}(U_1, U_2, U_3)$ *on* $(0, 1)^3$ *such that all its univariate margins are* Uniform$(0, 1)$ *and the following holds almost surely (with respect to joint Lebesgue measure):*

$$(3.8) \quad U_1 + U_2 + U_3 = 3/2 \quad and \quad \Phi^{-1}(U_1) + \Phi^{-1}(U_2) + \Phi^{-1}(U_3) = 0.$$

For $k$ even, say, $k = 4$, one may be tempted to use two pairs of quantile coupled variates, namely, $U_2 = 1 - U_1$, $U_4 = 1 - U_3$, where $U_1$ and $U_3$ are i.i.d., to make the $k = 4$ version of (3.8) hold. Consequently, for any $f$ $\rho_4^{(f)} = \rho_2^{(f)}/3$, where $\rho_2^{(f)} = \text{Corr}(f(U_1), f(1 - U_1))$ and $\rho_4^{(f)}$ is the correlation between any pair of $\{f(U_i), i = 1, \ldots, 4\}$. This implies $S_4^{(f)} = 1 + 3\rho_4^{(f)} = 1 + \rho_2^{(f)} = S_2^{(f)}$, and thus it is just a disguised version of quantile coupling with $k = 2$, not a real generalization to $k = 4$.

For lack of a universal strategy, we seek methods that are effective in common practice. In the exchangeable setting, the quest for EA is the same as that for minimizing the variance of the mean (and sum) $\bar{X}^{(k)} = (X_1 + \cdots + X_k)/k$ because [see (2.3)]

$$(3.9) \quad \rho_k \equiv \text{Corr}(X_i, X_j) = \frac{1}{k-1}\left[\frac{\text{Var}_k(\bar{X}^{(k)})}{\text{Var}(X_1)} - 1\right] \geq -\frac{1}{k-1},$$

where the subscript $k$ in both $\rho_k$ and $\text{Var}_k(\bar{X}^{(k)})$ emphasizes the dependence on the joint distribution $F^{(k)}(x_1, \ldots, x_k)$. In particular, if an $F^{(k)}$ is constructed such that $\bar{X}^{(k)} = constant$ (almost surely), then EA is achieved [hence (3.8)]. However, although this approach works for some common distributions such as uniform and normal, it is not always possible because the minimal value of $\text{Var}_k(\bar{X}^{(k)})$ may not achieve zero even when $k = 2$. For example, the minimal $\rho_2$ is $1 - \pi^2/6 = -0.645$ when both $X_1$ and $X_2$ are exponentially distributed with mean 1, as reported in Moran (1967).

For specific families of distributions, several methods may be available to achieve EA. Snijders (1984) explores various approaches for binary random variables. For a unimodal symmetric and differentiable density $p$ on $\mathbf{R}$, Rüschendorf and Uckelmann (2000) propose the following. Suppose the center of symmetry is zero and that $p_Q(x) = -xp'(x)$ is also a Lebesgue density on $\mathbf{R}$; let $Q \sim p_Q$. Then $X = QU \sim p$ for any $U \sim$ Uniform$(-1, 1)$ that is independent of $Q$. Consequently, for any set of $\{U_1, \ldots, U_k\}$, independent of $Q$, such that $U_i \sim$ Uniform$(-1, 1)$ and $\sum_{i=1}^{k} U_i = 0$, $\{X_1 = QU_1, \ldots, X_k = QU_k\}$ achieves EA with respect to $p$ because $\sum_{i=1}^{k} X_i = 0$. In fact, $\text{Corr}(X_i, X_j) = \text{Corr}(U_i, U_j) = -(k-1)^{-1}$ for any $i \neq j$. This construction, however, does not guarantee NA in general. As an alternative, we can draw i.i.d. $\{Q_1, \ldots, Q_k\}$ from $p_Q$, and then use $X_i = Q_i U_i$, $\forall i \in \{1, \ldots, k\}$. This will guarantee the NA property, but it sacrifices the EA property, because now we have $\text{Corr}(X_i, X_j) = -\frac{1}{k-1}[1 + CV^2(Q)]^{-1}$, where $CV(Q)$ is the coefficient of variation of $Q$. Nevertheless, when $CV(Q)$ is small, the loss of EA may be negligible for practical purposes. In general, if one has to make a choice between NA and EA, we recommend choosing NA, for it is preserved by all monotone transformations. With NA in place, the variance reduction factor is guaranteed to be at most 1 when the monotonicity assumption holds.

**4. Generating antithetic uniform variates.** Since generating uniform variates is, explicitly or implicitly, the most basic component of almost any simulation method, in this section we compare several methods for generating $k$-antithetically coupled uniform variates. For each method, we investigate whether it leads to NA and/or EA variables, and propose remedies whenever possible if it does not. We emphasize that although there are many different ways of achieving NA and/or EA [e.g., Bondesson (1983) and Gerow and Holbrook (1996)], no method dominates when $k \geq 3$, as demonstrated by Proposition 5. Any of these methods can be more suitable than others in a particular application. Nevertheless, the three methods described below are more or less representative of what has been proposed in the literature.

4.1. *The permuted displacement method.* This is a modified version of the one documented in Arvidsen and Johnsson (1982), which first generates an $r_1 \sim \text{Uniform}(0, 1)$, and then constructs

$$(4.1) \quad r_i = \{2^{i-2} r_1 + \tfrac{1}{2}\}, \qquad i = 2, \ldots, k-1, \quad \text{and} \quad r_k = 1 - \{2^{k-2} r_1\},$$

where $\{x\}$ is the fractional part of $x$. We find that the binary representation of this method makes it a bit easier to show that $\sum_{i=1}^{k} r_i = k/2$. Specifically, let $r_1 = (a_1, a_2, \ldots, a_m, \ldots)$ denote the (nonterminating) dyadic expansion of $r_1$, where $a_i = 0$ or $1$, that is, $r_1 = \sum_{i=1}^{\infty} a_i/2^i$. Then

$$r_2 = (1 - a_1, a_2, a_3, \ldots, a_m, \ldots),$$

$$r_3 = (1 - a_2, a_3, \ldots, a_{m+1}, \ldots),$$

$$\vdots$$

$$r_{k-1} = (1 - a_{k-2}, a_{k-1}, \ldots, a_{m+k-3}, \ldots),$$

$$r_k = (1 - a_{k-1}, 1 - a_k, \ldots, 1 - a_{m+k-2}, \ldots).$$

Therefore, the method creates negative correlation by displacing the binary digits of $r_1$. That all $r_i$'s are uniform is a direct consequence of a well-known result of Borel (1924), as discussed in detail by Billingsley (1986). The dyadic expansion representation shows clearly that the $r$'s are not exchangeable. A consequence of that is that they are *not* NA when $k \geq 3$. To prove this, we only need to show that when $k \geq 3$, $\{r_1, r_2\}$ are not NQD, and thus by Proposition 4 they cannot be NA. To see this, we note that for $k \geq 3$, $r_2 = r_1 + 0.5$ when $r_1 < 0.5$ and $r_2 = r_1 - 0.5$ when $r_1 \geq 0.5$. Consequently, for $0 \leq s < 0.5$ and $s + 0.5 < t \leq 1$, we have

$$P(r_1 \leq t, r_2 \leq s) = P(0.5 \leq r_1 \leq t, r_1 - 0.5 \leq s) = \min\{t, s + 0.5\} - 0.5 = s,$$

which is larger than $ts = P(r_1 \leq t)P(r_2 \leq s)$, and therefore (3.2) is violated.

However, it is easy to fix the nonexchangeability by using the simple permutation method. That is, let $S_k$ be the set of all possible permutations of $k$ objects. Pick a random $\sigma \in S_k$ and define $U_i = r_{\sigma(i)}$. Clearly $\sum_i U_i = \sum_i r_i = k/2$ and thus $\mathrm{Corr}(U_i, U_j) = -(k-1)^{-1}$ for any $i \neq j$. Furthermore, it can be shown that:

THEOREM 3. *For $k = 3$, $\{U_1, U_2, U_3\}$ constructed by the permuted displacement method are PNA.*

The proof is given in Craiu and Meng (2002) and is omitted both because of space limitation and because the approach used was rather brute force. Indeed, we are unable to prove or disprove Theorem 3 for $k \geq 4$. Nevertheless, the result indicates that the exchangeability can play an important role in achieving NA/PNA. An astute reader might wonder how permuting indexes can be helpful since in MCMC our estimates typically are sample averages, which are invariant to the independent permutations of the sample indexes. However, one must keep in mind that in our use of the antithetic variates, the $U$'s are not used just once in the end, but throughout the whole generated sequence and the final estimates are not invariant to the permutations of all the indexes that occur along the sequence. This is perhaps easiest to see for CFTP, as each draw can depend, in principle, on an arbitrarily large number of independent copies of the $k$-tuple of $U$'s.

4.2. *Multidimensional normal method.* A common way to manipulate correlations, especially in the engineering literature, is through the multivariate normal distribution. For our purposes, we can first generate $(Z_1, \ldots, Z_{k-1})^\top \sim N(\mathbf{0}, \Sigma)$ where $\Sigma_{ij} = -(k-1)^{-1}$ if $i \neq j$ and $\Sigma_{ii} = 1$, and then let $Z_k = -(Z_1 + Z_2 + \cdots + Z_{k-1})$. Finally, if we are interested in uniform deviates, we can use $U_i = \Phi(Z_i)$ for all $i \in \{1, \ldots, k\}$, where $\Phi$ is the CDF of $N(0, 1)$. The random variables $\{U_1, \ldots, U_k\}$ are NA following Joag-Dev and Proschan (1983), who proved that multivariate normals with nonpositive pairwise correlations are NA. The result then follows immediately from Proposition 3 because $U_i$ is a monotone function of $Z_i$.

This method, however, does not achieve EA. Although $\text{Corr}(Z_i, Z_j) = -(k-1)^{-1}$, the nonlinear transformation $U_i = \Phi(Z_i)$ causes an increase in correlation. Specifically, for any $\{Z_1, Z_2\}$ bivariate normal with correlation $\rho$, the corresponding correlation between $\Phi(Z_1)$ and $\Phi(Z_2)$ is given by

$$(4.2) \qquad \frac{\int \Phi(\rho t)\Phi(t)\phi(t)\,dt - 1/4}{1/12} = \frac{3}{\pi}\rho + \frac{2\sqrt{\pi}-3}{16\pi}\rho^3 + O(\rho^5)$$
$$= 0.955\rho + 0.01\rho^3 + O(\rho^5),$$

where $\phi$ is the standard normal density. Consequently, $\text{Corr}(U_1, U_2)$ is larger than the minimum possible value $-(k-1)^{-1}$, though the loss of EA may not be important for many practical applications in view of (4.2). Equation (4.2) is a quantitative version of a result of Lancaster (1957), who proved that for any functions $g_1$ and $g_2$, $|\text{Corr}(g_1(Z_1), g_2(Z_2))| \leq |\rho|$ as long as the left-hand side is finite, and the equality holds if and only if $g_i(z) = z$, $i = 1, 2$ (almost surely). Lancaster's result also implies Proposition 5 if we require $\{Z_1, Z_2, Z_3\}$ to be jointly normal. Nevertheless, Proposition 5 is a stronger result, as it shows that it is impossible to find any such trivariate distribution, not just trivariate normal, to simultaneously preserve EA as in (3.8).

An issue with real impact in computation is the requirement for a highly reliable and efficient subroutine to evaluate the function $\Phi$. Otherwise, the use of a not-highly accurate approximation becomes problematic in large replications with arbitrarily many arguments because once in a while $|Z|$ can be too large for $\Phi(Z)$ to be evaluated appropriately. For that reason, we did not use this method to simulate uniform deviates in the simulated examples in Section 2. However, we used it to generate antithetic truncated normal deviates in the probit example presented there.

The normal method is also of interest because it is one that many will likely attempt, especially for generating antithetic normal deviates in high-dimensional settings. An example where such implementation results in acceleration of a classical MCMC method is provided by Craiu (2004) in the context of Multiple-Try Metropolis with antithetic proposals.

4.3. *Iterative Latin hypercube sampling.* The method described in Section 4.1 achieves EA but whether it achieves NA is an open question, and the method given in Section 4.2 achieves NA but not EA (when used for generating the uniform variates). To achieve both, we propose to use *iterative Latin hypercube sampling* (ILHS), which is an iterative version of the Latin hypercube sampling (LHS), a well-known scheme in quasi Monte Carlo; see McKay, Beckman and Conover (1979), Stein (1987), Owen (1992), Loh (1996), Iman (1999) and Helton and Davis (2003), among others.

For any given $k \geq 2$, our iterative procedure can be described by the following steps:

1. Set $\mathcal{U}_0^k = (U_0^{(1)}, \ldots, U_0^{(k)})^\top$, where $\{U_0^{(j)}\}_{1 \le j \le k}$ are i.i.d. Uniform$(0,1)$.
2. For $t = 0, 1, 2, \ldots$, let $\mathcal{K}_t = (\sigma_t(0), \ldots, \sigma_t(k-1))^\top$ be a permutation of $\{0, 1, \ldots, k-1\}$, independent of all previous draws, and let

$$(4.3) \qquad \mathcal{U}_{t+1}^k = \frac{1}{k}(\mathcal{K}_t + \mathcal{U}_t^k),$$

where $\mathcal{U}_t^k = (U_t^{(1)}, \ldots, U_t^{(k)})^\top$ for $t \ge 0$.

The case of $t = 1$ corresponds to the original LHS. For general $t$ we have the following result.

THEOREM 4.  *For any $t \ge 0$, $i, j \in \{1, \ldots, k\}$, $i \ne j$, we have:*

  (i) $U_t^{(i)} \sim \text{Uniform}(0,1)$.
  (ii) $\text{Corr}(U_t^{(i)}, U_t^{(j)}) = -\frac{1}{k-1}(1 - \frac{1}{k^{2t}})$.
  (iii) $\{U_{t_1}^{(1)}, \ldots, U_{t_k}^{(k)}\}$ *are NA for any finite* $\{t_1, \ldots, t_k\}$.

PROOF.  (i) For $t = 0$, the result obviously holds. Suppose $U_t^{(1)} \sim \text{Uniform}(0,1)$; then

$$P(U_{t+1}^{(1)} \le s) = P(U_t^{(1)} < ks - \sigma_t(1)) = \frac{1}{k} \sum_{j=0}^{[ks]} P(U < ks - j) = s$$

for any $s \in (0,1)$. The result thus holds by induction.

(ii) Let $S_t = \sum_j U_t^{(j)} = \mathbf{1}^\top \mathcal{U}_t^k$. Then from the recursion formula (4.3) we have $\text{Var}(S_t) = \text{Var}(S_{t-i})/k^{i+1}$ for all $1 \le i \le t$, which implies, by (i) and the exchangeability of $\{U_t^{(1)}, \ldots, U_t^{(k)}\}$,

$$k + k(k-1)\,\text{Corr}(U_t^{(i)}, U_t^{(j)}) = \frac{k}{k^{2t}},$$

which is just (ii).

(iii) Because any permutation distribution is NA, (4.3) defines an antithetically coupled joint Markov chain with $k$ marginal Markov chains. Furthermore, the marginal updating function from $U_{t-1}^{(j)}$ to $U_t^{(j)}$ is monotone. The result then follows directly from part (ii) of Theorem 1 because $\{U_0^{(1)}, \ldots, U_0^{(k)}\}$ are i.i.d. and thus NA.                    □

For most practical purposes, Theorem 4 is all we need as it proves that for any $T$ we choose, $\{U_T^{(1)}, \ldots, U_T^{(k)}\}$ are NA. Furthermore, as long as $T$ is not too small (e.g., $T \ge 5$), they practically achieve EA because the relative loss of EA is $k^{-2T}$, which approaches zero very fast (it is less than 0.02% even for $k = 3$ and $T = 5$). However, in theory $\{U_T^{(1)}, \ldots, U_T^{(k)}\}$ achieves EA only for $T = \infty$. This requires showing first that $\{U_\infty^{(1)}, \ldots, U_\infty^{(k)}\}$ are well defined and second that they are NA. The first question is easy to answer but the second is not, mainly because the support of $\{U_\infty^{(1)}, \ldots, U_\infty^{(k)}\}$ is a "Cantor dust" type of fractal, as investigated in Craiu and Meng ([2002]).

In practice, one must choose $T$ and $k$. We have already seen that large $T$ results in more extreme antithesis, while the simulations performed in Section [2] show that as $k$ increases, the efficiency increases too. The following result, proved in Appendix [A.2], further shows that, *without* taking into account the computational cost, the larger $T$ the better.

THEOREM 5. *For any monotone $h \in L^2[0,1]$, the correlation* $\mathrm{Corr}(h(U_{T+1}^{(1)}), h(U_{T+1}^{(2)}))$ *is decreasing as a function of $T$. That is, for any $T \geq 0$,*

$$(4.4) \qquad \mathrm{Corr}(h(U_{T+1}^{(1)}), h(U_{T+1}^{(2)})) \leq \mathrm{Corr}(h(U_T^{(1)}), h(U_T^{(2)})).$$

However, $T$ and $k$ cannot be increased indefinitely in practice. Intuitively, when $k$ becomes large, the choice of $T$ should become less important because a large $k$ means we have a deep stratification even with $T = 1$. Consequently, further stratification within each stratum, which is essentially what each new iteration does, becomes less important. The following result, proved in Appendix [A.3], provides a theoretical support of this intuition by showing that for large $k$, the impact of $T$ is negligible.

THEOREM 6. *For any $h \in L^2[0,1]$ and a fixed $T$, $\{h(U_T^{(1)}), \ldots, h(U_T^{(k)})\}$ achieves EA asymptotically as $k \to +\infty$, that is,*

$$(4.5) \qquad \mathrm{Corr}(h(U_T^{(1)}), h(U_T^{(2)})) = -\frac{1}{k-1} + o(k^{-1}).$$

This result implies that asymptotically, as $k \to \infty$, any ILHS iteration achieves EA for any square-integrable function, not just monotone functions. In practice, however, $k$ must be finite, and often quite small for the sake of computational cost. For fixed $k$, the following result, proved in Appendix [A.4], shows that even with $k$ as small as 3, $T$ does not need to be large in order to achieve practically the same result as $T = \infty$, at least for all monotone estimand functions.

THEOREM 7. *Let $D_{h_1,h_2}(t, t+m)$ be the Kolmogorov–Smirnov distance between the two-way (marginal) joint CDF of $\{h_1(U_{t+m}^{(1)}), h_2(U_{t+m}^{(2)})\}$ and of $\{h_1(U_t^{(1)}), h_2(U_t^{(2)})\}$, where $h_l$, $l = 1, 2$, are nondecreasing functions. Then*

$$(4.6) \qquad D_{h_1,h_2}(t, t+m) \leq k^{-(t-1)}(k-1)^{-(t+2)}.$$

This implies that if we take

$$T \geq \frac{d - 2 \log_{10}(k-1) + \log_{10} k}{\log_{10}(k(k-1))},$$

then $D_{h_1,h_2}(T,\infty) \leq 10^{-d}$. In particular, as long as $T \geq 5$, $D_{h_1,h_2}(T,\infty) <$ 0.0001 for any $k \geq 3$. Thus, taking $T$ between 5 and 10 is enough for almost any practical purpose. The generality of Theorem 7 should be emphasized since the same bound in (4.6) holds for any monotone $h_1$ and $h_2$.

4.4. *Variance reduction factors for indicator functions.* As an attempt to find some theoretical support of the empirical findings reported in Section 2 that the cost-effective choice of $k$ appears to be somewhere between 3 and 10, we report here our findings of the theoretical value of $S_k^{(f)}$ as a function of $k$ when $f = \mathbb{1}_{\{x \leq c\}}$. We choose the indicator function both because of its analytical tractability and its practical relevance (e.g., for estimating CDF), and because it serves as a building block for general functions. We start with the LHS method, for which

$$(4.7) \qquad S_k^{(f)}(c) = 1 + (k-1)\frac{F(c,c) - c^2}{c(1-c)} = \frac{(1 - \{kc\})\{kc\}}{kc(1-c)},$$

where $\{a\}$ denotes the fractional part of $a$. This result follows from the expression of the joint CDF, $F(c,c) = P(U_1 \leq c, U_2 \leq c)$, given by (A.6) from Appendix A.2; note (A.5) implies that the same $S_k^{(f)}(c)$ holds for any ILHS.

The left panel of Figure 8, realized using the freely available software RGL developed by Duncan Murdoch, plots $S_k^{(f)}(c)$ as a function of both $k$ (up to 30) and $c \in (0,1)$. Its fascinating shape reveals that as long as $c$ is not too close to 0 or 1, $S_k^{(f)}(c)$ will be rather small. This is more clearly seen in the first two rows of Figure 9, where $S_k^{(f)}(c)$ is plotted against $c$ for given $k$. It is intuitive that when $c$ approaches 0 or 1, the effect of antithetic coupling fades because the indicator function approaches the constant function. This can be formalized by considering that for $1/k \leq c \leq (k-1)/k$ the maximum of $S_k^{(f)}(c)$ over this range, as shown in Appendix A.5, is given by

$$(4.8) \qquad S_k^* = \frac{k}{3k - 4 + 2\sqrt{2(k-1)(k-2)}},$$

corresponding to the dashed lines in the first two rows in Figure 9. The use of a large $k$ is seen to be twice advantageous: first, it decreases the $S_k^*$ and second, it shrinks the area (of $c$) where the antithetic coupling is not effective. The plots also show clearly that $k = 2$ is least effective.

The $S_k^*$ of (4.8) decrease from $S_3^* = 1/3$ to $S_\infty^* = (3 + 2\sqrt{2})^{-1} \approx 0.172$. However, the use of large $k$ also increases the computational cost. Striving for a balance, we consider $R_k = (S_3^* - S_k^*)/(S_3^* - S_\infty^*)$, a measure of the relative gain in efficiency obtained when we increase from $k = 3$ to a finite $k > 3$ instead of $k = \infty$. Table 2 gives, for a few $\alpha \in [0.5, 1]$, $k_\alpha = \min\{k : R_k \geq \alpha\}$. It is seen that with $k = 9$ we already have achieved 80% of the additional gain. The rapid decrease of $S_k^{(f)}(c)$ as a function of $k$, for fixed $c$, can also be seen in the last two rows of Figure 9.

TABLE 2
*The minimum $k_\alpha$ for achieving $\alpha\%$ of maximal possible gain
in efficiency over using $k = 3$*

| $\alpha$ | 50% | 60% | 70% | 80% | 90% | 95% | 99% |
|---|---|---|---|---|---|---|---|
| $k_\alpha$ | 5 | 6 | 7 | 9 | 17 | 32 | 152 |

The above exercise can be repeated for the multivariate normal method, that is, when $(X_1, X_2, \ldots, X_k)$ is a multivariate normal vector with $\mathrm{Corr}(X_i, X_j) = -(k-1)^{-1}$ and $N(0, 1)$ margins. The calculation of $S_k^{(f)}$ is a bit more involved than in the LHS case. Specifically, with the orthogonal transformation $Z = X_1 + X_2$, $W = X_1 - X_2$ (i.e., $Z$ and $W$ are independent), we obtain

$$\Phi_k(c, c) \equiv P(X_1 \leq c, X_2 \leq c) = P(Z - 2c \leq W \leq -Z + 2c, Z \leq 2c)$$

$$(4.9) \qquad = 2 \int_{-\infty}^{2c/\sqrt{2(1+\rho_k)}} \Phi\left(\frac{2c - z\sqrt{2(1+\rho_k)}}{\sqrt{2(1-\rho_k)}}\right) \phi(z)\, dz$$

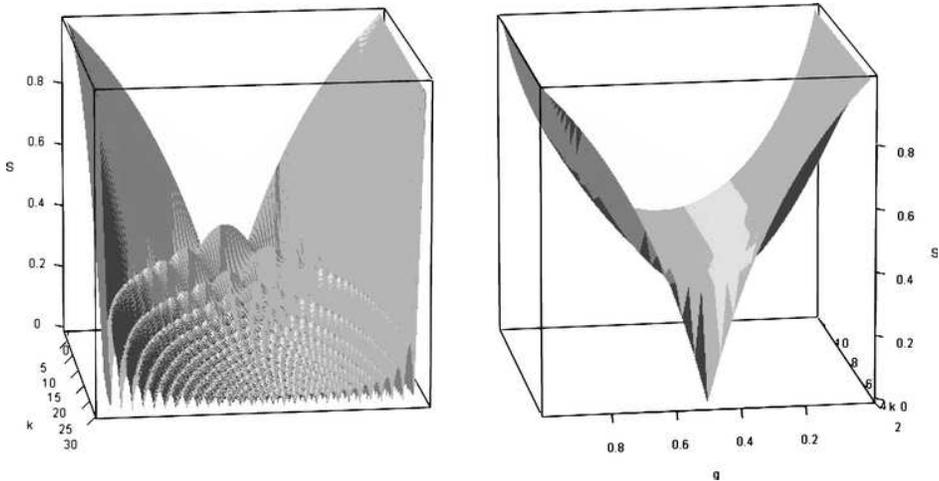

FIG. 8. *Variance reduction for indicator functions for* Uniform$(0, 1)$ *random variables generated using LHS* (left panel) *and for normal variates* (right panel). *Note for better visualization the left panel is seen from back, with the c axis hidden.*

$$- \Phi\left(\frac{2c}{\sqrt{2(1+\rho_k)}}\right),$$

where $\rho_k = -(k-1)^{-1}$. Consequently, we can compute $S_k^{(f)}(c)$ via

$$(4.10) \qquad S_k^{(f)}(c) = 1 + (k-1)\frac{\Phi_k(c,c) - \Phi^2(c)}{\Phi(c)(1-\Phi(c))}.$$

In the normal case, we expect a behavior of $S_k^{(f)}(c)$ similar to that in the uniform case when $c$ is close to the limits of the range, $(-\infty, \infty)$. However, for plotting purposes we use a one-to-one transformation of $c$, $g = \Phi(c)$, and we plot $S_k^{(f)}(c)$ against $g$; this is the same as the VRF for the uniform variates via the normal approach, namely, if $U_i = \Phi(X_i)$. One can notice in the right panel of Figure 8 that everywhere but in a region of $c$ around the distribution's center, $S_k^{(f)}(c)$ is decreasing in $k$. However, the decrease is much less abrupt than in the uniform case, which indicates that it will take a larger $k$ to reach the same relative efficiency as given in Table 2. This suggests that generalizing the findings in Table tb:eff1 to other situations is by no means automatic.

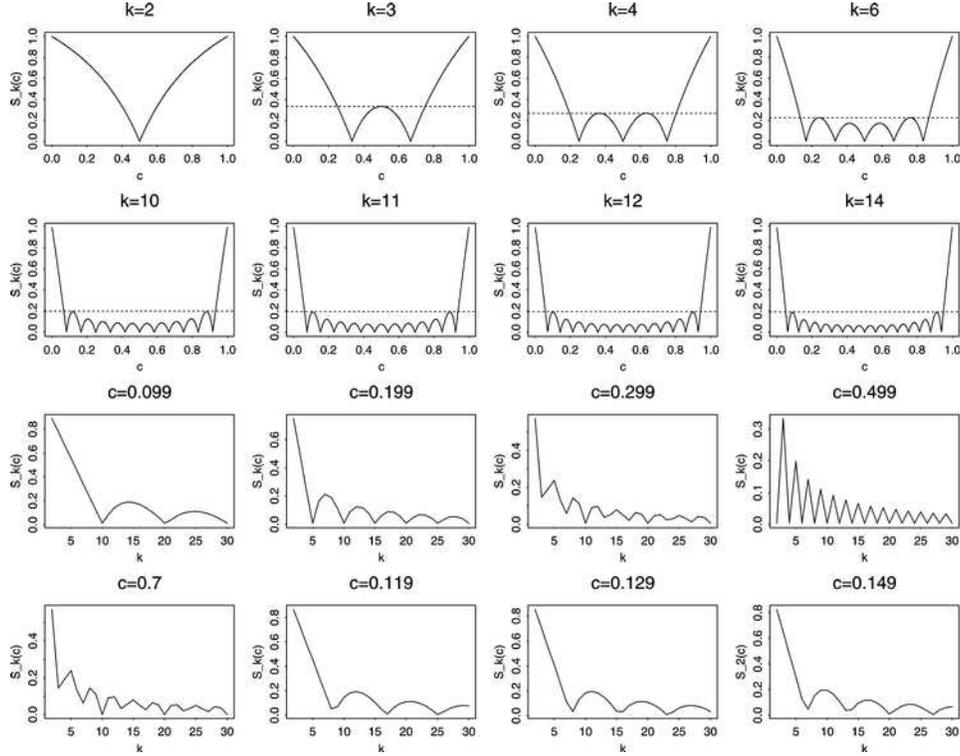

FIG. 9. *Variance reduction factor for indicator function under LHS method, as a function of c for different values of k* (first two rows), *and as a function of k for different values of c* (last two rows).

This example also shows that it is not always true that the larger $k$, the better. For example, when $c = 0$, using $k = 2$ completely eliminates the variance because of the symmetry in $f = \mathbb{1}_{\{x \leq 0\}}$; see the right panel of Figure 8. This can also be verified directly from (4.10), which becomes

$$S_k^{(f)}(0) = 1 + (k-1)\left[8 \int_0^\infty \Phi\left(z\sqrt{1 - \frac{2}{k}}\right)\phi(z)\,dz - 3\right],$$

which is zero when $k = 2$. However, as $k \to \infty$ we can show that $S_k^{(f)}(0) \to 1 - 2/\pi = 0.3634$. This reminds us that although there are good rules of thumb for general practice, in terms of choosing $k$ as well as other choices (e.g., the generating methods) presented in this paper, one should not adopt them blindly without examining the special structures of the problem at hand.

# APPENDIX

A.1. *Proof for the unit circle example.* We first note that (3.4) implies, for any $t \geq 1$,

$$(A.1) \qquad (\theta_t^X + \theta_t^Y) \bmod (2\pi) = -(\theta_{t-1}^X + \theta_{t-1}^Y) \bmod (2\pi).$$

From the right panel of Figure 7, $\widehat{MX_0} = \widehat{QN} = \widehat{QY_0} + \widehat{Y_0N}$, where $\widehat{AB}$ denotes the counter-clockwise arc between points $A$ and $B$ on the unit circle. Furthermore, $\widehat{MX_0} = \theta_0^X - \pi$, $\widehat{QY_0} = 2\tau$ and $\widehat{Y_0N} = 2\pi - \theta_0^Y$. Combining these identities, we obtain $\theta_0^X + \theta_0^Y = 2\tau + 3\pi$. Consequently, (A.1) implies (3.5) by noting that $\theta_t^X + \theta_t^Y \in (0, 4\pi)$ and the "alternating" nature of (A.1), which is also clear from Figure 7, and thus the four possible lines in (3.5) are reachable from each other by $(X_t, Y_t)$. Similar arguments apply to other possible initial configurations of $X_0$ and $Y_0$ (e.g., when $X_0$ and $Y_0$ are on different sides of the horizontal axes).

To prove $\rho_t(\tau) = -\cos(2\tau)$, we observe that (3.5) implies $\sin^2(\theta_t^X + \theta_t^Y) = \sin^2(2\tau)$ for any $t \geq 0$. Using the identity $\sin^2(\alpha + \beta) = \cos^2(\alpha) + \sin^2(\beta) - 2\cos(\alpha)\cos(\beta)\cos(\alpha + \beta)$, we then obtain

$$(A.2) \qquad \cos^2(\theta_t^Y) + \cos^2(\theta_t^X) + 2\cos(\theta_t^Y)\cos(\theta_t^X)\cos(2\tau) = \sin^2(2\tau).$$

Thus the orbit of $(X_t^{(1)}, Y_t^{(1)}) \equiv (\cos(\theta_t^X), \cos(\theta_t^Y))$ is an ellipse. Taking expectations on both sides of (A.2) and using $\mathrm{E}[\cos^2(\theta_t^Y)] = \mathrm{E}[\cos^2(\theta_t^X)] = 1/2$, we obtain $2\mathrm{E}[\cos(\theta_t^Y)\cos(\theta_t^X)]\cos(2\tau) = -\cos^2(2\tau)$. Together with $\mathrm{Var}(\cos^2(\theta_t^X)) = \mathrm{Var}(\cos^2(\theta_t^Y)) = 1/2$, this yields $\rho_t(\tau) = -\cos(2\tau)$ when $\tau \neq \pi/4$. For $\tau = \pi/4$, (3.5) gives $\theta_t^Y = (2k+1)\pi \pm \pi/2 - \theta_t^X$, and thus $\cos(\theta_t^Y) = \pm\sin(\theta_t^X)$. Consequently, $2\mathrm{E}[\cos(\theta_t^Y)\cos(\theta_t^X)] = \pm\mathrm{E}[\sin(2\theta_t^X)] = 0$, so $\rho_t(\tau) = -\cos(2\tau)$ still holds.

A.2. *Proof of Theorem* 5. Since the (marginal) distribution of $U_T^{(j)}$ is Uniform$(0, 1)$ and thus does not depend on $T$, by the Hoeffding identity (1.1), inequality (4.4) becomes immediate if we can establish

$$(\text{A}.3) \qquad F_{T+1}^{(h)}(u, v) \le F_T^{(h)}(u, v) \qquad \forall\, (u, v),$$

where $F_t^{(h)}(u, v)$ is the joint CDF of $(h(U_t^{(1)}), h(U_t^{(2)}))$.

We first show that (A.3) holds for $h(x) = x$ by mathematical induction. Given the exchangeability, we can assume $u \le v$, and thus $[ku] \le [kv]$, where $[x]$ denotes the integer part of $x$. Let $p_{ij} = P(U_t^{(1)} \le ku - i, U_t^{(2)} \le kv - j)$; then by the recursion (4.3),

$$
\begin{aligned}
(\text{A}.4) \quad F_{t+1}(u, v) &= P(U_t^{(1)} \le ku - \sigma_t^{(1)}, U_t^{(2)} \le kv - \sigma_t^{(2)}) \\
&= \frac{1}{k(k-1)} \sum_{i=0}^{[ku]} \sum_{j=0}^{[kv]} p_{ij} \mathbb{1}_{\{i \ne j\}}.
\end{aligned}
$$

To evaluate this expression, we consider the following possibilities according to the value of $(i, j)$:

(A) When $i \le [ku] - 1$ and $j \le [kv] - 1$, $p_{ij} = P(U_t^{(1)} \le 1, U_t^{(2)} \le 1) = 1$; there are $[ku]([kv] - 1)^+$ such pairs of $(i, j)$'s within $0 \le i \ne j \le k$. Recall $(x)^+ = \max\{x, 0\}$.

(B) When $i \le [ku] - 1$ and $j = [kv]$, $p_{ij} = P(U_t^{(2)} \le kv - [kv]) = \{kv\}$, where $\{x\} = x - [x]$; there are $[ku]$ such pairs of $(i, j)$'s.

(C) When $i = [ku]$ and $j \le [kv] - 1$, $p_{ij} = P(U_t^{(1)} \le ku - [ku]) = \{ku\}$; there are $[kv]$ such pairs of $(i, j)$'s when $[ku] = [kv]$, but only $[kv] - 1$ such pairs when $[ku] < [kv]$ because of the requirement that $i \ne j$.

(D) When $i = [ku]$ and $j = [kv]$, $p_{ij} = P(U_t^{(1)} \le ku - [ku], U_t^{(2)} \le kv - [kv]) = F_t(\{ku\}, \{kv\})$. There is no such pair when $[ku] = [kv]$ and one such pair when $[ku] < [kv]$.

Putting these four possibilities together, we have

$$
\begin{aligned}
&F_{t+1}(u, v) \\
(\text{A}.5) \quad &= \begin{cases} 0, & \text{if } [kv] = 0, \\ \dfrac{ku(kv - 1) - \{ku\}(\{kv\} - 1)}{k(k-1)}, & \text{if } 0 < [kv] = [ku], \\ \dfrac{ku(kv - 1) - \{ku\}\{kv\} + F_t(\{ku\}, \{kv\})}{k(k-1)}, & \text{if } [kv] > [ku]. \end{cases}
\end{aligned}
$$

From (A.5), if (A.3) holds for $T = 0$ the rest follows immediately by induction. For $t = 1$, since $F_0(\{ku\}, \{kv\}) = \{ku\}\{kv\}$, we have from (A.5)

$$(A.6) \quad F_1(u, v) = \begin{cases} 0, & \text{if } [kv] = 0, \\ \dfrac{ku(kv-1) - \{ku\}(\{kv\}-1)}{k(k-1)}, & \text{if } 0 < [kv] = [ku], \\ \dfrac{u(kv-1)}{k-1}, & \text{if } [kv] > [ku]. \end{cases}$$

Among the three expressions above, the second one is the largest as a function of $u$ and $v$. It is easy to check that this second function is less than or equal to $F_0(u, v) = uv$ if and only if

$$(A.7) \quad \{ku\}(1 - \{kv\}) \le ku(1 - v).$$

Using $[ku] = [kv]$, (A.7) is also equivalent to $\{kv\}(u - \{ku\}) \le [ku](1-u)$, which is obvious when $u \le \{ku\}$. When $0 \le \{ku\} < u$, $\{ku\}(1 - \{kv\}) \le u(1 - \{kv\})$. But then $u(1 - \{kv\}) \le ku(1-v)$ is equivalent to $[kv] \le k - 1$, which is obviously true except when $v = 1$. But the $v = 1$ case is trivial because then $u = 1$ in order to maintain $[ku] = [kv] = k$, and thus $\{ku\} = 0$, and hence (A.7) holds for all $0 \le u \le v \le 1$. This proves $F_1(u, v) \le F_0(u, v)$ for all $(u, v)$ and hence (A.3) by induction.

To show that (A.3) holds for any nondecreasing $h$, let $x_w = \sup\{x : h(x) \le w\}$ for any given $w$. Then $\{x : h(x) \le w\} = A_w^{(h)}(x)$, where $A_w^{(h)}(x) = \{x : x \le x_w\}$ if $h(x_w) \le w$ and $A_w^{(h)}(x) = \{x : x < x_w\}$ if $h(x_w) > w$. This means that the probability calculations of event $\{U : h(U) \le w\}$ are the same as those for either $\{U : U \le x_w\}$ or $\{U : U < x_w\}$. This allows us to go from $F_T^{(h)}(u, v)$ to $F_T(x_u, x_v)$, for which we already have proved the desired result. A technical complication here is that, depending on the continuity properties of $h$ at $x_u$ and $x_v$, one or two "$\le$" operations in the definition of $F_T(x_u, x_v) = P(U_T^{(1)} \le x_u, U_T^{(2)} \le x_v)$ may need to be replaced by "$<$". However, as far as (A.3) is concerned, such modifications are immaterial because they do not affect any part of (A)–(D). Or mathematically, we have, for any given $(u, v)$,

$$\begin{aligned} F_{T+1}^{(h)}(u, v) &= P(h(U_{T+1}^{(1)}) \le u, h(U_{T+1}^{(2)}) \le v) \\ &= E[A_u^{(h)}(U_{T+1}^{(1)}) \cap A_v^{(h)}(U_{T+1}^{(2)})] \le E[A_u^{(h)}(U_T^{(1)}) \cap A_v^{(h)}(U_T^{(2)})] \\ &= P(h(U_T^{(1)}) \le u, h(U_T^{(2)}) \le v) = F_T^{(h)}(u, v), \end{aligned}$$

where the middle inequality holds because (A.3) holds with $h(x) = x$, possibly with the aforementioned modifications from "$\le$" to "$<$" in the definition of the bivariate CDF.

A.3. *Proof of Theorem* 6. For any $h \in L^2[0, 1]$, Stein (1987) showed that Theorem 6 is true for the original LHS, that is, when $T = 1$. Therefore, for any monotone $h \in L^2[0, 1]$, Theorem 6 is a consequence of Theorem 5, because the smallest possible value of $\mathrm{Corr}(h(U_T^{(1)}), h(U_T^{(2)}))$ is $-1/(k-1)$, as seen in (3.9). More specifically, for any $T > 1$, any possible limit of $-(k-1)\mathrm{Corr}(h(U_T^{(1)}), h(U_T^{(2)}))$, as $k \to \infty$, must be bounded below and above by 1, and hence can only be 1. For nonmonotone $h \in L^2[0, 1]$, the proof turns out to be much more technical. It essentially requires going through Stein's (1987) original arguments but using the more complicated bivariate CDF via (A.5); details are given in Craiu (2001).

A.4. *Proof of Theorem* 7. First we prove (4.6) when both $h_1$ and $h_2$ are identity functions. By (A.5),

$$(A.8) \quad D(t, t+1) = \sup_{[ku] \neq [kv]} \left\{ \left| \frac{F_t(\{ku\}, \{kv\})}{k(k-1)} - \frac{F_{t-1}(\{ku\}, \{kv\})}{k(k-1)} \right| \right\}$$

$$(A.9) \quad \leq \sup_{u,v} \left\{ \left| \frac{F_t(u, v)}{k(k-1)} - \frac{F_{t-1}(u, v)}{k(k-1)} \right| \right\} = \frac{D(t-1, t)}{k(k-1)}.$$

This allows us to conclude that

$$(A.10) \qquad\qquad D(t, t+1) \leq \frac{D(0, 1)}{[k(k-1)]^t},$$

where $D(0, 1) = \sup_{u,v} \{ |F^{(1)}(u, v) - F^{(0)}(u, v)| \} = \sup_{u,v} \{ uv - F^{(1)}(u, v) \}$. By (A.6), $D(0, 1)$ is the maximum value of the following three suprema:

(i) $\displaystyle \sup_{[ku]=[kv]=0} uv = \frac{1}{k^2}$;

(ii) $\displaystyle \sup_{0 < [kv]=[ku]} \left\{ uv - \frac{ku(kv-1) - \{ku\}(\{kv\}-1)}{k(k-1)} \right\}$

   $\displaystyle = \sup_{0 < [kv]=[ku]} \left\{ \frac{ku(1-v) - \{ku\}(1-\{kv\})}{k(k-1)} \right\} \leq \frac{1}{k-1}$;

(iii) $\displaystyle \sup_{[kv] \neq [ku]} \left\{ uv - \frac{ku(kv-1)}{k(k-1)} \right\} = \sup_{[kv] \neq [ku]} \left\{ \frac{u(1-v)}{k-1} \right\} = \frac{1}{k-1}$.

Putting all these facts together, we have

$$D(t, t+m) \leq \sum_{i=t}^{t+m-1} D(i, i+1)$$

$$\leq \sum_{i=t}^{t+m-1} \frac{1}{k^i(k-1)^{i+1}} = \frac{1}{k^t(k-1)^{t+1}} \frac{1 - [k(k-1)]^{-m}}{1 - [k(k-1)]^{-1}},$$

which is less than $k^{-(t-1)}(k-1)^{-(t+2)}$ when $k \geq 2$. The extension to nondecreasing functions $h_1, h_2$ is immediate by defining $h^{-1}(u) = \sup\{s : h(s) \leq u\}$ and proceeding in the same fashion as in Appendix A.2. That is, for any nondecreasing functions $h_1$ and $h_2$,

$$D_{h_1, h_2}(t, t+m) \leq \sup_{u,v} |F_{t+m}(h_1^{-1}(u), h_2^{-1}(v)) - F_t(h_1^{-1}(u), h_2^{-1}(v))|$$

$$\leq D(t, t+m),$$

where the definition of $F_{t+m}(x_u, x_v)$ and $F_t(x_u, x_v)$ may need to be modified as in Appendix A.2.

A.5. *Maximum reduction factor in the uniform case.* To maximize (4.7) subject to $1 \leq [kc] \leq k-1$, let $i_c = [kc]$ and $f_c = \{kc\}$. By symmetry, the maximum of (4.7) occurs at either $i_c = 1$ or $i_c = k-2$. When $i_c = 1$, (4.7) can be expressed as

$$(A.11) \qquad S_k^{(f)}(f_c) = \frac{k(1-f_c)f_c}{(1+f_c)(k-1-f_c)}.$$

Straightforward differentiation then shows that the maximizer must satisfy $(k-3)f_c^2 + 2(k-1)f_c - (k-1) = 0$, which only has one acceptable solution $f_c = \frac{\sqrt{k-1}}{\sqrt{2k-4}+\sqrt{k-1}}$. Therefore the maximizer of (4.7) is $c_1^* = \frac{1}{k}(1 + \frac{\sqrt{k-1}}{\sqrt{2k-4}+\sqrt{k-1}})$ with the corresponding maximal value of $S_k^{(f)}$ given in (4.8). The maximum from $i_c = k-2$ is the same, with the maximizer $c_2^* = \frac{1}{k}(k-2 + \frac{\sqrt{k-1}}{\sqrt{2k-4}+\sqrt{k-1}})$.

**Acknowledgments.** We thank Vanja Dukic, Steve Lalley, Duncan Murdoch, Art Owen, Michael Stein, Stephen Stigler, Lei Sun and Mike Wichura for helpful exchanges. We also thank the reviewers for a set of exceptionally thorough and helpful comments.

Department of Statistics
University of Toronto
100 St. George Street
Toronto, Ontario
Canada M5S 3G3
e-mail: craiu@ustat.toronto.edu

Department of Statistics
Harvard University
1 Oxford Street
Cambridge, Massachusetts 02138
USA
e-mail: meng@stat.harvard.edu